\documentclass[11pt,letterpaper]{amsart}
\usepackage{amssymb, amscd, charter,epic, eepic}
\numberwithin{equation}{section}

\def\AA{{\mathbb A}}
\def\CC{{\mathbb C}}  
\def\GG{{\mathbb G}} 

\def\HH{{\mathbb H}}

\def\PP{{\mathbb P}}
 
\def\RR{{\mathbb R}}
\def\TT{{\mathbb T}} 
\def\UU{{\mathbb U}}  
\def\VV{{\mathbb V}} 
\def\ZZ{{\mathbb Z}}

\def\odd{{\rm odd}} 
\def\even{{\rm even}} 
\def\spec{{\rm sp}} 
\def\orb{{\rm orb}} 

\def\lin{{\rm f}}

\def\reg{{\rm reg}}

\def\tor{{\rm tor}} 

\def\sw{{SW}}
\def\wbar{\overline{W}}
\def\swbar{{\overline{SW}}}
\def\sar{{S\! Art}}
\def\Ar{Art}
\def\sarbar{\overline{\sar}}

\def\arbar{\overline{Art}}

\def\Ione{\text{I}_1}
\def\Itwo{\text{I}_1}

\def\G{\Gamma} 

\def\bs{\backslash}

\def\Ccal{{\mathcal C}}

\def\Ecal{{\mathcal E}} 
 
\def\Hcal{{\mathcal H}} 
\def\Kcal{{\mathcal K}}
\def\Lcal{{\mathcal L}}
\def\Mcal{{\mathcal M}}
\def\Scal{{\mathcal S}}

\def\Mtilde{\widetilde{M}}
\def\Mcalvec{\Mcal_{3,\vec{1}}}
\def\Gvec{\Gamma_{3,\vec{1}}}

\def\Scalhat{{\hat \Scal}}

\def\Ocal{{\mathcal O}}
\def\Scal{{\mathcal S}}

\def\Vcal{{\mathcal V}}  
\def\Xcal{{\mathcal X}}
 
\def\half{\frac{1}{2}}

\newcommand\aut{\operatorname{Aut}}
\newcommand\bl{\operatorname{Bl}}

\newcommand\GL{\operatorname{GL}}
\newcommand\SL{\operatorname{SL}}

\newcommand\PGL{\operatorname{PGL}}

\newcommand\Hom{\operatorname{Hom}}

\newcommand\pic{\operatorname{Pic}}
\newcommand\sym{\operatorname{Sym}}
\newcommand\ord{\operatorname{ord}}  

\newtheorem{theorem}{Theorem}[section]
\newtheorem{lemma}[theorem]{Lemma}
\newtheorem{proposition}[theorem]{Proposition}
\newtheorem{corollary}[theorem]{Corollary}

\theoremstyle{definition}
\newtheorem{definition}[theorem]{Definition}

\theoremstyle{remark} 
\newtheorem{remark}[theorem]{Remark}

\begin{document}

\title{Artin groups and the fundamental groups of some moduli
spaces} 

\author{Eduard Looijenga}

\address{Faculteit Wiskunde en Informatica, University of
Utrecht, P.O.\ Box 80.010, NL-3508 TA Utrecht, The
Netherlands}

\email{looijeng@math.ruu.nl}
\subjclass[2000]{14J10, 20F36, 20F34, 14H10}

\begin{abstract}
We define for every affine Coxeter graph a 
certain factor group of the associated Artin group
and prove that some of these groups appear as orbifold
fundamental groups of moduli spaces. Examples are the 
moduli space of nonsingular cubic algebraic surfaces and the 
universal nonhyperelliptic smooth genus three curve.
We use this to obtain a simple presentation of 
the mapping class group of a compact genus three 
topological surface with connected
boundary. This leads to a modification of Wajnryb's presentation
of the mapping class groups in the higher genus case that
can be understood in algebro-geometric terms.
\end{abstract}

\thanks{Support by the Institut des Hautes \'Etudes Scientifiques
and the Mittag-Leffler Institute is gratefully acknowledged}
\maketitle

\section*{Introduction}

\noindent
If $G$ is a group (having the discrete topology)
acting on a connected space $X$ and $p\in X$ is a base point,
then the set of pairs $(\alpha ,g)$ with $g\in G$ and $\alpha$
a homotopy class in $X$ of paths from $p$ to $gp$ forms a group
under the composition law $(\alpha ,g).(\beta ,h)=
(\alpha (g_*\beta), gh)$ (we adopt the geometric convention for
the multiplication law for fundamental groupoids:
if $\gamma, \delta$ are paths in $X$ such that
$\delta$ begins where $\gamma$ ends, then 
$\gamma\delta$ stands for $\gamma$ followed by $\delta$). 
This is the {\it equivariant
fundamental group}  $\pi_1^G (X,p)$ of $(X,p,G)$.
The projection on the second  component defines a surjective 
homomorphism $\pi_1^G (X,p)\to G$ with kernel $\pi_1(X,p)$. 

Following an insight of Borel, one can also define $\pi_1^G (X,p)$ as the fundamental group of
an actual space $X_G$: choose a universal $G$-bundle $EG\to BG$ 
(with $G$ acting on the right of $EG$) and take $X_G:=EG\times^G X$. 
The advantage of this definition is that it makes also sense if 
$G$ is a topological group 
acting continuously on $X$. Since $X_G\to BG$ is a 
fiber bundle with fiber $X$, the Serre sequence gives rise 
to an exact sequence
\[
\pi_1(G,e)\to \pi_1(X,p)\to \pi_1^G (X,p)\to \pi_0(G)\to 1 
\]
with  the first map induced by $g\in G\mapsto gp\in X$
(here we used that the Serre sequence of $EG\to BG$ yields an isomorphism
$\pi_n(BG,*)\cong\pi_{n-1}(G,e)$). In all cases that we consider here,
$X$ is a manifold and $G$ is a Lie group acting properly on $X$ with finite
stabilizers. It is then convenient to think of $X_G$ as the orbit space 
$G\bs X$ with a bit of extra structure given at the orbits 
with nontrivial stabilizer, as formalized by the notion of 
an orbifold. This is reflected by a
change in terminology and notation: we call $\pi_1^G (X,p)$ 
the \emph{orbifold fundamental group} of $X_G$ and write it also
as $\pi_1^\orb (X_G,p)$.  The kernel $K$ of the $G$-action on $X$ naturally 
appears as a normal subgroup of $\pi_1^\orb (X_G,p)$ and the
corresponding quotient can be identified with $\pi_1^\orb (X_{G/K},p)$.
We will call this the \emph{reduced orbifold fundamental group} of $X_G$.

A fundamental example taken from algebraic geometry is the following: for
integers $n\ge 0$ and $d\ge 3$, let $H^n(d)$ denote the dicriminant complement
in $\PP (\sym^d\CC^{n+2})$, that is, the Zariski open subset
of $\PP (\sym^d\CC^{n+2})$ that parameterizes the {\it nonsingular} 
hypersurfaces in $\PP^{n+1}$ of degree $d$. It supports the universal
degree $d$ hypersurface $X^n(d)\subset \PP^{n+1}_{H^n(d)}$ and on both
these algebraic manifolds the group $\PGL(n+2)$ acts properly 
(and hence with finite stabilizers, as the group is affine).
So the orbifold fundamental group of $H^n(d)_{\PGL(n+2)}$ is
a quotient of the fundamental group of $H^n(d)$ by 
the image of the fundamental group of $\PGL(n+2)$ (which may be 
identified with the center $\mu_{n+2}$ of its universal cover $\SL(n+2)$)
and similarly for $X^n(d)_{\PGL(n+2)}$. 
For $n\ge 3$ a nonsingular hypersurface in $\PP^n$ of degree $d\ge 3$
is simply connected and has in general no nontrivial automorphism.
This implies that for such $n$ and $d$, the projection 
$X^n(d)_{\PGL(n+2)}\to H^n(d)_{\PGL(n+2)}$ induces an isomorphism
of orbifold fundamental groups. So the separate consideration
of $X^n(d)_{\PGL(n+2)}$ is only worthwhile for 
$n=0$ ($d$-element subsets of $\PP^1$) or $n=1$ (nonsingular plane curves
of degree $d$).  

Here is what we believe is known about these  groups.
There is the classical case of the coarse moduli space 
$d$-element subsets on a smooth rational curve, $H^0(d)_{\PGL (2)}$, 
whose orbifold fundamental group is isomorphic to the
mapping class group of the $2$-sphere less $d$ points
(where these points may get permuted) and for which
a presentation is due to Birman \cite{birman} (although this 
presentation does not treat the $d$ points on equal footing). 
Its covering $X^0(d)_{\PGL (2)}$ can be interpreted as
the coarse moduli space of  $(d-1)$-element subsets on an affine line, 
and its orbifold fundamental group is in fact the usual Artin braid group
with $d-1$ strands modulo its center. Another case with a classical flavor is
the moduli space of smooth genus one curves endowed with a divisor class of degree $3$, $H^1(3)_{\PGL (3)}$.   Dolgachev and Libgober showed \cite{dollib} that the orbifold 
fundamental group of $H^1(3)_{\PGL (3)}$ is the semidirect product
$\SL(2,\ZZ )\ltimes (\ZZ/(3))^2$. It is then easy to see that the orbifold 
fundamental group of the universal such curve, $X^1(3)_{\PGL (3)}$,
is the semidirect product $\SL(2,\ZZ )\ltimes \ZZ^2$
with the projection $X^1(3)_{\PGL(3)}\to H^1(3)_{\PGL (3)}$ inducing the obvious
map between these semidirect products.
Finally, a  concrete presentation for the fundamental group of
any $H^n(d)$ has been recently obtained by Michael L\"onne \cite{loenne}.

In this paper we shall, among other things, find a presentation for the orbifold fundamental groups of the moduli spaces of Del Pezzo surfaces of degree 3, 4 and 5 (in the last case there are no moduli and the orbifold fundamental group is simply the automorphism group of such a surface). We should perhaps recall that for $d=3$ this is also the orbifold fundamental group of the moduli space of cubic surfaces $H^2(3)_{\PGL(4)}$. In all these cases the presentatation is obtained as a quotient of  some Artin group of affine type (in this case $\hat E_6$, $\hat D_5$, $\hat A_4$). These quotients can in fact be defined in a uniform manner for all such Artin groups. For $\hat E_7$ (degree 2) it yields the orbifold fundamental group of the universal quartic curve $X^1(4)_{\PGL(3)}$ and we also have a result  for $\hat E_8$ (degree 1).

By pushing our methods a bit further we also
get an efficient  presentation of the mapping class group 
of a compact genus three surface with a single
boundary component relative to its boundary. This 
is obtained by purely algebro-geometric means. The presentation 
that we find is new and is similar (but not identical) to a 
M.~Matsumoto's simplification of the 
presentation due to Hatcher-Thurston, Harer and Wajnryb.
It leads also to a geometrically meaningful  presentation
of the mapping class group of any surface with a single boundary 
component and of genus at least three.

My interest in these questions were originally aroused by conversations
with Domingo Toledo, while both of us were staying at the Institut des
Hautes \'Etudes Scientifiques (IH\'ES) in Bures.
In fact, it was his asking for a presentation
of the orbifold fundamental group of the moduli space of 
nonsingular cubic surfaces that is at the origin of this paper.
He and his collaborators Allcock and Carlson have obtained
a beautiful arithmetic representation of this fundamental group
\cite{allcock}. 

The essential part of this work was done in while
I was enjoying the hospitality of the IH\'ES in the Spring of 1997.
The present paper is a substantial reworking of an earlier version,
which was carried out during a visit of the Mittag-Leffler
institute on 2007. I am very grateful to both institutions for their 
support and hospitality and for offering the ambiance so 
favourable for research. 

\smallskip
{\it Plan of the paper.} Section $2$, 
after a brief review of some basic facts about Artin groups, 
introduces certain factor groups of an Artin group of 
affine type. This leads up to the
formulation and discussion of our principal results.
In Section $3$ we investigate these factor groups
in more detail. Here the discussion is still
in the context of group theory, but in 
Section $4$ we introduce natural orbifolds associated
to Artin groups of affine type having such a group as orbifold 
fundamental group. In Section $5$ we occupy ourselves with the following
question: given a family of arithmetic genus one curves
whose general fiber is an irreducible rational
curve with a node and a finite set of sections, what happens
to the relative positions of the images of these sections
in the Picard group of the fibers, when we approach 
a special fiber? The answer is used to construct 
over the orbifolds introduced in Section $4$ a family of what we have baptized
Del Pezzo triples (of a fixed degree). When the degree is
at least three, this moduli space can be understood as
the normalization of the dual of the universal Del Pezzo
surface of that degree. It contains an 
orbifold as constructed in Section $4$ as the complement
of a subvariety of codimension two and so has the same
orbifold fundamental group. This leads to a proof of 
our main results with the exception of those which concern the 
mapping class group of genus three. The last section  
is devoted to that case.

\tableofcontents

\smallskip
\textit{Some notational conventions.}
If a reflection group $W$ acts properly on a complex manifold $M$, 
then $M^\circ$ will denote the locus of $x\in M$ not fixed
by a reflection in  $W$. The use if this notation is of course 
only permitted if there is no ambiguity on $W$.

If $C$ is a complete curve, then $\pic_k(C)$ stands for the set of divisor classes
of degree $k$. We also follow the custom to denote the identity component of a topological group $G$ by  $G^0$. So $\pic_0(C)\subset \pic (C)^0$ and 
this is an equality when $C$ is irreducible.

Finally, we write $\Delta$ for the spectrum of $\CC [[t]]$
and we denote its closed point by $o$ and its generic point
by ${\Delta^*}$. (The symbol $\Delta$ with a subscript such as
$\Delta_\G$, stands for something entirely 
different---namely an element of a 
certain group---but no confusion should arise.)

\section{Statement of the results}\label{sect:results}

\subsection*{Artin groups and their reductions}
Let $\G$ be a graph with vertex set $I$. Assume that 
$\G$ has no loops and that any two vertices are connected by at
most countably many edges. Then $\G$ defines a presentation
of group $\Ar_\G$ with generators $\{t_i : i\in I \}$
indexed by the vertex set  and relations given as follows:
if $i,j\in I$ are  distinct and connected by $k_{ij}$
edges, then impose the relation 
\[ 
t_it_jt_it_j\cdots = t_jt_it_jt_i\cdots 
\] 
with on both sides $k_{ij}+2$ letters,
whenever that number is finite (no relation is imposed when
$k_{ij}=\infty$). It is called the {\it Artin group} attached 
to $\G$.\index{Artin group} 
Such groups were first investigated by 
Brieskorn-Saito \cite{briesaito} and Deligne \cite{deligne}.
If we are also given its set of generators
$\{t_i\}_{i\in I}$, then we shall shall refer to the
pair $(\Ar_\G, \{t_i\}_{i\in I})$ as an {\it Artin system}\index{Artin system}.
The {\it length homomorphism} \index{length homomorphism}
$\ell :\Ar_\G\to \ZZ$ is the homomorphism 
which takes on each generator $t_i$ the value 1. 
If $\G'\subset\G$ is a full subgraph of $\G$, then according 
to Van der Lek \cite{lek:thesis}, \cite{lek:paper} the obvious
homomorphism $\Ar_{\G'}\to \Ar_\G$ is injective and so we may
regard $\Ar_{\G'}$ as a subgroup of $\Ar_\G$. 

Imposing also the relations $t_i^2=1$ (all $i$) gives
the {\it Coxeter group} $W_\G$ defined by $\G$. If $s_i$ denotes
the image of $t_i$ in $W_\G$, then $(W_\G, \{s_i\}_{i\in
I})$ is a Coxeter system in the sense of Bourbaki
\cite{bourb}. It is shown there that $W_\G$ has a natural faithful
representation in $\RR^I$ with $s_i$ acting as a reflection and
and $s_is_j$ ($i\not= j$) acting as a transformation of 
order $k_{ij}+2$. So $\G$ can be
reconstructed from the Coxeter system and hence also from
the Artin system. The automorphism group of $\G$ acts
on $\Ar_\G$ and we can therefore form the semidirect product
$\Ar_\G\rtimes \aut (\G )$.

We say that the graph $\G$ is {\it of finite type} 
\index{finite type!graph of}
if $W_\G$ is finite. 
Brieskorn-Saito and
Deligne have shown that for connected $\G$ of finite type 
the group of $u\in\Ar_\G$ with the  property that 
conjugation by $u$ preserves the generating
set $\{ t_i\}_{i\in I}$ is infinite cyclic. The homomorphism $\ell$
is nontrivial on this subgroup, so there is  a unique 
generator $\Delta_\G$ on which
$\ell$ takes a positive value. We call it the {\it Garside element}.
\index{Garside element} Its square is always central. (The geometric 
meaning of these facts will be recalled in \ref{finitecoxeter}.)
So the inner automorphism defined by the Garside element is in fact
an automorphism of the  system (of order
at most two) and defines an 
automorphism of $\G$ of the same order, the 
{\it canonical involution}\index{involution!canonical} 
of $\G$. The canonical involution is the identity if
the associated Coxeter group contains minus the identity in
its natural representation and in the other cases ($A_{k\ge 2}$, 
$D_\odd$, $E_6$, $I_\odd$), it is the  unique involution of $\G$.
If $\G$ is of finite type, but not necessarily connected, then
$\Delta_\G$ is by definition the product of the Garside elements
of the connected components and the canonical involution
of $\Delta_\G$ is induced by conjugation with that element
(and so the canonical involution on each component). 

We say that a connected graph $\G$ is of {\it affine type}
\index{affine!type, connected} if
the natural representation of $W_\G$ in $\RR^I$ 
leaves invariant an affine hyperplane in the dual of $\RR^I$
on which it acts faithfully and properly. (A characterization in
terms of group theory is that $W_\G$  
contains a nontrivial free abelian subgroup of finite index
and is not a product of two nontrivial subgroups.)
The classification of such graphs can be found in \cite{bourb},
Ch.~VI, \S\ 4, Thm.\ 4.

If $\G$ is connected, of finite type 
and given as the underlying graph of a connected
Dynkin diagram (this excludes the types $H_3$, $H_4$ and $\Itwo(m)$
for $m=5,7,8,\dots $ and distinguishes the cases $B_l$ and $C_l$), 
then there is a 
natural {\it affine completion}\index{affine!completion} 
$\G\subset \hat\G$, where $\hat\G$ is of affine type,
has one vertex more than $\G$, and contains $\G$ as a full subgraph.
All connected graphs of affine type so arise (though not
always in a unique manner, see below).
The construction is canonical and so the natural involution of 
$\G$ extends to $\hat\G$ (fixing the added vertex), but for what 
follows it is important to point out that this involution is 
{\it not} induced by conjugation with $\Delta_\G$. \\

\emph{In the remainder of this section we assume that  $\G$ 
is connected and of affine type}. The groups associated 
to $\G$ that we are about 
to define appear in the statements of our main results. This
is the sole reason for introducing them now; a 
fuller discussion will be given in Section \ref{artin:affine}.

Given any $i\in I$, then the full subgraph 
$\G_i$ on $I-\{ i\}$ is always of finite type and so we have defined
a Garside element $\Delta_{\G_i}$ which we here simply denote by $\Delta_i$.
We say that $i$ is {\it special}\index{special!vertex}
if $\G$ is the natural completion of $\G_i$. As we just noticed, 
a special vertex defines an involution $g_i$ of $\G$. We shall write $I_\spec$ 
for the set of special vertices. It turns out that $I_\spec$ 
is an orbit of the automorphism group of $\G$ and so the isomorphism 
type of $\G_i$ is determined by that of $\G$. If $\G$ is of type
$\hat A_l$ (i.e., a polygon with $l+1$ vertices), then all vertices are special. 
But for the other connected graphs of affine type there 
are only few of these: the graphs in question are trees and their special 
vertices are ends (i.e., lie on a single edge). To be precise,  
for those of type $\hat D_l$ we get all $4$ ends, 
for $\hat B_l$ and $\hat C_l$ $2$ ends, for $\hat E_6$, $\hat E_7$, $\hat E_8$ 
the ends on the longest branches (resp.\ $3$, $2$, $1$ in number) 
and in the remaining cases $\hat F_4$ and $\hat G_2$ only one. 

There is a similar property for pairs 
$(i,j)$ of distinct special vertices in $I$:
then there is an involution $g_{ij}$ of $\G$ which interchanges
$i$ and $j$ and is the canonical involution on the 
full subgraph $\G_{ij}$ of $\G$ spanned by $I-\{ i,j\}$.
If $\G$ is not a polygon, than $\G_{ij}$ is obtained from $\G$ by 
removing two ends and hence is still connected.
The elements $g_i$ and $g_{ij}$ generate a subgroup 
$S(\G )\subset\aut(\G)$, which is in fact all of $\aut(\G)$ unless $\G$ is of
type $\hat D_\even$. 

\begin{definition} 
We define the \emph{reduced Artin group} $\arbar_\G$ as
the quotient of the semidirect product $\Ar_\G\rtimes S(\G )$
by the imposing the  relations: 
\begin{align*}
\Delta_i g_i \equiv 1\quad&\text{for $i$ special.} \tag{$1$}\\
\Delta_{ij}g_{ij}\equiv 1\quad&\text{for $i,j$ distinct and special.}\tag{$2_m$}\\
\end{align*}
Here $\Delta_{ij}$ is the Garside element for the
Artin subgroup defined by the full subgraph on $I-\{ i,j\}$.
\end{definition} 

A few observations are in order. It is clear from the definition that
the evident homomorphism $\Ar_\G\to\arbar_\G$ is surjective. A description
of $\arbar_\G$ as a quotient of $\Ar_\G$ amounts to saying that 
conjugation with any word $w_\Delta$ in the  Garside elements $\Delta_i$ and 
$\Delta_{ij}$ permutes the generating set $\{ t_i\}_{i\in I}$ as 
the corresponding word $w_g$ in the $g_i$'s and the $g_{ij}$'s and that 
$w_\Delta$ represents the identity element precisely when $w_g$ does.
In particular, 
\begin{align*}
\Delta_i t_i \equiv t_i\Delta_i,\quad \Delta_i^{\ord g_i}\equiv 1
\quad &\text{for $i$ special,} \\
\Delta_{ij} t_i \equiv t_j\Delta_{ij},\quad \Delta_{ij}^{\ord g_{ij}}\equiv 1
\quad  & \text{for $i,j$ distinct and special}.
\end{align*}
For a pair $(i,j)$ as above, the identity
$t_i\equiv \Delta_{\G_{ij}}t_j\Delta_{\G_{ij}}^{-1}$ shows that 
the image of $t_i$ in $\arbar_\G$ can be expressed in terms of the 
images of the other generators $t_j$, $j\not= i$, so that
$\arbar_\G$ is in fact a quotient of $\Ar_{\G_i}$. But this 
breaks the symmetry and as a consequence the defining relations are 
no longer easy to describe. 

As already mentioned, the affine Coxeter group $W_\G$ 
is naturally realized as an affine transformation group. 
The quotient of this affine transformation group by its translation 
subgroup is a finite (Coxeter) group $W_\G^\lin$ acting 
naturally in a vector space $V$. The center $Z(W_\G^\lin)$ of 
$W_\G^\lin$ is the intersection of $W_\G^\lin$ with 
$\{ \pm 1_{V}\}$.  We shall show that the quotient of 
$\sarbar_\G$ obtained by imposing the relations 
$t_i^2\equiv 1$ for all $i\in I$ can be identified with 
$W_\G^\lin/Z(W_\G^\lin)$, reason for us to denote the latter
group simply by $\wbar_\G$. We will see that in general the 
homomorphism $\Ar_{\G}\to W_\G^\lin$ has a
nontrivial kernel. An exception is the case $\hat A_l$:

\begin{proposition}\label{prop:smallgroup}
For $l>1$, the natural homomorphism
$\overline{Ar}_{\hat A_l}\to \overline{W}_{\hat A_l}\cong\Scal_{l+1}$
is an isomorphism.
\end{proposition}

We will prove this (relatively easy) assertion in Section \ref{artin:affine}.

\subsection*{Some orbifold fundamental groups}
We can now state some of our main results. These concern the moduli spaces
of Del Pezzo surfaces of degree $d\le 5$. Recall that such a surface can be obtained
by blowing up $9-d$ points of the projective plane in general position (but is not
naturally given this way). For $d=2,3,4$ the moduli space in question is the space of projective equivalence classes of respectively smooth quartic curves in $\PP^2$, of cubic surfaces in $\PP^3$ and of complete intersections of two quadrics in $\PP^4$.
For $d=5$ there are no moduli since the Del Pezzo surfaces of degree 5 are mutually isomorphic. But that does not imply that its  orbifold fundamental group is trivial:
this will be the (finite) automorphism group of that surface.

\begin{theorem}\label{thm:main} 
The orbifold fundamental group of the moduli space of Del Pezzo surfaces of degree 
resp.\ 5,4 and 3 is isomorphic to a reduced Artin group of type resp.\ $\hat A_4$, $\hat D_5$  and $\hat E_6$. The natural surjection of this group on the finite Weyl group of the corresponding type (resp.\ $A_4$, $D_5$ and $E_6$) describes the representation of orbifold fundamental group on the 
integral homology of the Del Pezzo surface. 
\end{theorem}

Notice that Proposition \ref{prop:smallgroup} implies that the automorphism group of a degree 5 Del Pezzo surface acts faithfully (via the finite Weyl group of type $A_4$)
on its homology.

As remarked above, we can change the presentation of 
the reduced Artin group of type $\hat E_6$ into one which does
not involve one of the special generators. We thus
recover a theorem of Libgober \cite{libgo}, which says
that the orbifold fundamental group of the moduli space of cubic surfaces
is a quotient of an Artin group of type $E_6$. Our result says even which
quotient, but as we noted, such a description does not do
justice to the $\hat E_6$-symmetry.

For Del Pezzo moduli in degree 1 and 2 the situation is a little different.
The graph of type $\hat E_7$ has two special vertices. Let us write simply
$E_7$ and $E'_7$ for the corresponding subdiagrams and $E_6$ for 
their intersection (this labeling is of course of a selfreferencing 
nature).

\begin{theorem}\label{thm:main2} 
The orbifold fundamental group of the universal smooth quartic plane curve,  
given up to projective equivalence is isomorphic to a reduced Artin group of type 
$\hat E_7$: it is the quotient of the Artin group of type $\hat E_7$ by the relations 
\begin{enumerate}
\item[(i)] $\Delta_{E_7}\equiv 1$,
\item[(ii)] $t_1\Delta_{E_6}\equiv\Delta_{E_6} t_7$ and
$\Delta_{E_6}^2\equiv 1$.
\end{enumerate}
The natural surjection of this group on the the Weyl group of type $E_7$ 
modulo its center describes the representation of orbifold fundamental group 
on the homology of a smooth quartic curve 
with $\ZZ /2$-coefficients.
\end{theorem}

Perhaps we should point out that relation (ii) in this theorem implies
that  $\Delta_{E'_7}=\Delta_{E_6}\Delta_{E_7}\Delta_{E_6}^{-1}$ so 
that  $\Delta_{E'_7}\equiv 1$ by relation (i) and we indeed
obtain the reduced Artin group of type $\hat E_7$.

This theorem leads to a new simple presentation of the mapping
class group of genus three. Recall that any smooth
(complex) projective curve of genus three is either a
quartic curve or hyperelliptic. The hyperelliptic curves
form in the moduli space of smooth projective genus three
curves a divisor, so we expect the orbifold fundamental
group of the universal smooth complex genus three curve to
be a quotient of $\sarbar_{\hat E_7}$. This is the case: if 
$A_7\subset \hat E_7$ denotes the unique subgraph of that type, then
we have:

\begin{theorem}\label{thm:genus3}
The orbifold fundamental group of the universal smooth complex genus three curve 
is naturally the quotient of the reduced affine Artin group of type $\hat E_7$  by the relation $\Delta_{A_7}\equiv\Delta_{E_6}$. 
\end{theorem}

\begin{remark}
We also find something of interest in the case $d=1$: the reduced Artin group
of type $\hat E_8$ is the orbifold fundamental group of
pairs $(S,K)$, where  $S$ is a Del Pezzo surface of degree 1 and $K$ is
a singular member of its anticanonical pencil (there are in general 12 such
members).
\end{remark}

\subsection*{Presentation of the mapping class group}
Take a compact oriented  surface $\Sigma_g$ of genus $g$
with connected boundary $\partial \Sigma_g$ (so  $\partial \Sigma_g$
is diffeomorphic to a circle). We can consider the group of orientation preserving
diffeomorphisms of $\Sigma_g$ and its subgroup of  diffeomorphisms that
are the identity on the boundary. We follow Harer in denoting the groups of connected components by $\G_g^1$ and $\G_{g,1}$ respectively (although in view of a convention in algebraic geometry, a notation like $\G_{g,1}$ resp.\ $\G_{g,\vec{1}}$ would perhaps have been more logical).
The evident map $\G_{g,1}\to \G_g^1$ is surjective and its kernel is the (infinite cyclic) group generated by the Dehn twist along the boundary. Their Eilenberg-Maclane spaces have algebro-geometric incarnations that we now recall.

If $g\ge 1$, then the moduli space   of pairs $(X,p)$, where $X$ is a 
complete nonsingular genus $g$ curve and $p\in X$ is an orbifold 
(Deligne-Mumford stack) $\Mcal_{g,1}$. Its orbifold fundamental group
is $\G_g^1$ and its orbifold universal covering is contractible. 
If in addition is given a nonzero tangent vector to $X$ at $p$, then
the corresponding moduli  space, which we propose to denote by  $\Mcal_{g,\vec{1}}$, is fine and smooth: this is because an automorphism of a projective nonsingular curve of positive genus which fixes a nonzero tangent vector must be the identity. Notice that it comes with an action of  $\CC^\times$  (acting as scalar multiplication on the tangent vector). That action is proper (but not free: finite nontrivial  isotropy groups do occur) and its orbit space may be identified with $\Mcal_{g,1}$. The smooth variety 
$\Mcal_{g,\vec{1}}$ is an Eilenberg-MacLane space for  $\G_{g,1}$: it has the latter
as  its fundamental group and its universal covering is contractible. The projection $\Mcal_{g,\vec{1}}\to \Mcal_{g,1}$ induces the natural map
$\G_{g,1}\to \G_g^1$ and the positive generator of a fiber $(\cong \CC^\times$) corresponds to a Dehn twist along the boundary.
The fundamental group of  $\Mcal_{3,\vec{1}}$  can also be presented as a quotient of $\Ar_{\hat E_7}$:

\begin{theorem}\label{mappingclass} 
The fundamental group of $\Mcalvec$ (which is also is the mapping class group $\Gvec$ of a compact oriented genus three surface with a circle as boundary)
is isomorphic to the quotient of
$\Ar_{\hat E_7}$ by the following relations.
\begin{enumerate} 
\item[(i)] $\Delta_{E_7}\equiv\Delta_{E_6}^2$, 
\item[(ii)] $\Delta_{A_7}\equiv\Delta_{E_6}$. 
\end{enumerate} 
\end{theorem}

\begin{figure}
\setlength{\unitlength}{0.00056868in}%
\begingroup\makeatletter\ifx\SetFigFont\undefined%
\gdef\SetFigFont#1#2#3#4#5{%
  \reset@font\fontsize{#1}{#2pt}%
  \fontfamily{#3}\fontseries{#4}\fontshape{#5}%
  \selectfont}%
\fi\endgroup%
{\renewcommand{\dashlinestretch}{30}
\begin{picture}(8535,3835)(0,-10)
\put(7530.000,2440.500){\arc{1984.597}{5.2878}{7.2786}}
\put(1567.500,2913.000){\arc{1530.662}{0.4242}{2.7174}}
\put(1567.500,2159.250){\arc{1238.318}{3.5727}{5.8520}}
\put(4222.500,2913.000){\arc{1530.662}{0.4242}{2.7174}}
\put(4222.500,2159.250){\arc{1238.318}{3.5727}{5.8520}}
\put(6697.500,2913.000){\arc{1530.662}{0.4242}{2.7174}}
\put(6697.500,2159.250){\arc{1238.318}{3.5727}{5.8520}}
\put(1005.000,2440.500){\arc{1984.597}{2.1462}{4.1370}}
\put(3073.977,3823.227){\arc{1100.456}{0.1832}{1.5689}}
\put(5416.022,3823.228){\arc{1100.455}{1.5727}{2.9584}}
\thicklines
\put(480.395,2433.395){\arc{961.036}{3.3931}{6.2215}}
\put(2882.239,2099.642){\arc{1716.906}{3.5215}{5.8461}}
\put(5428.500,2206.500){\arc{1503.242}{3.5543}{5.9349}}
\put(7831.641,2129.719){\arc{1570.904}{3.7107}{5.7807}}
\put(4222.500,1900.500){\arc{497.041}{4.8030}{7.7633}}
\thinlines
\put(4245,3723){\ellipse{1260}{180}}
\thicklines
\put(1590,2462){\ellipse{1620}{900}}
\put(4245,2463){\ellipse{1620}{900}}
\put(6674,2460){\ellipse{1620}{900}}
\put(465,978){\blacken\ellipse{128}{128}}
\put(465,978){\ellipse{128}{128}}
\put(4155,978){\blacken\ellipse{128}{128}}
\put(4155,978){\ellipse{128}{128}}
\put(4155,78){\blacken\ellipse{128}{128}}
\put(4155,78){\ellipse{128}{128}}
\put(7935,978){\blacken\ellipse{128}{128}}
\put(7935,978){\ellipse{128}{128}}
\put(2940,978){\blacken\ellipse{128}{128}}
\put(2940,978){\ellipse{128}{128}}
\put(1770,978){\blacken\ellipse{128}{128}}
\put(1770,978){\ellipse{128}{128}}
\put(6675,978){\blacken\ellipse{128}{128}}
\put(6675,978){\ellipse{128}{128}}
\put(5415,978){\blacken\ellipse{128}{128}}
\put(5415,978){\ellipse{128}{128}}
\thinlines
\path(465,3273)(3075,3273)
\path(465,3273)(3075,3273)
\path(465,3273)(3075,3273)
\path(465,3273)(3075,3273)
\path(465,3273)(3075,3273)
\path(465,3273)(3075,3273)
\path(465,3273)(3075,3273)
\path(465,3273)(3075,3273)
\path(465,3273)(3075,3273)
\path(465,3273)(3075,3273)
\path(5415,3273)(8025,3273)
\path(5415,3273)(8025,3273)
\path(420,1653)(8115,1653)(8115,1608)
\thicklines
\path(465,978)(7890,978)(7935,978)
\path(4110,978)(4155,978)(4155,78)
\end{picture}
}
\caption{The genus three surface with its $\hat E_7$-configuration}
\end{figure}
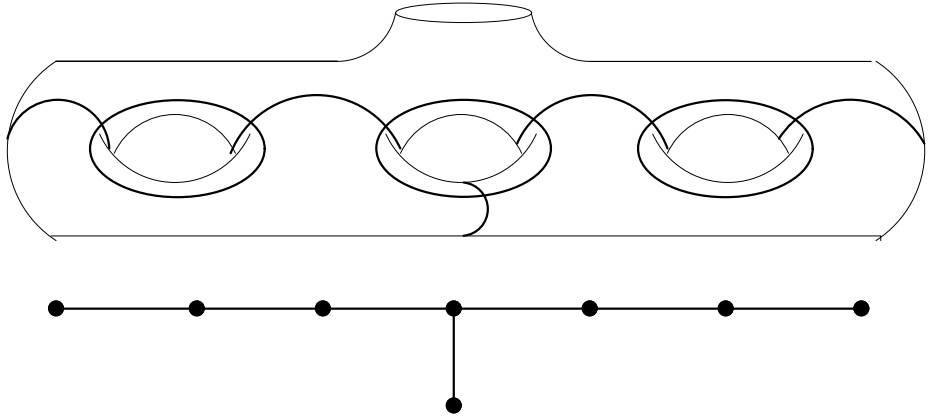

The proof will show that the Artin generators map to
Dehn twists along embedded circles which meet like a $\hat E_7$-diagram 
as in  Figure 1 (the Dehn twists along two curves 
$\alpha,\beta$ commute if the curves are disjoint, whereas
in case they meet simply in a single point their isotopy classes $D_\alpha, D_\beta$ 
satisfy Artin's braid relation $D_\alpha D_\beta D_\alpha=D_\beta D_\alpha D_\beta$).
It can be verified that the Garside elements $\Delta_{A_7}$ and
$\Delta_{E_6}$ themselves can be represented by the same 
homeomorphism; this homeomorphism is the composite of an involution
of $\Sigma_{3,1}$ with three (interior) fixed points (in Figure 1 this is the reflection with respect to the vertical axis of symmetry) and a
half Dehn twist along $\partial \Sigma_{3,1}$ to ensure that it
leaves $\partial \Sigma_{3,1}$ pointwise fixed. The Garside element $\Delta_{E_7}$
is the  Dehn twist along $\partial \Sigma_{3,1}$.

A rather explicit, but still somewhat involved presentation of
this group is due to Wajnryb \cite{wajn}. Makoto Matsumoto
\cite{mats} recently deduced with the help of a computer from that presentation a
leaner one which differs from ours by just one relation (he
also obtains the first relation) and involves Artin groups. Nevertheless, at present
the equivalence of his and our presentation has not yet
been established. In any case, the following geometric description of 
the generators and relations in the above theorem
should make a direct comparison with the one of Wajnryb and 
Matsumoto possible.  Another difference with his work is 
that we get the generators and relations almost handed
to us by some classical algebraic geometry so that the presentation that
we get is a very natural one.

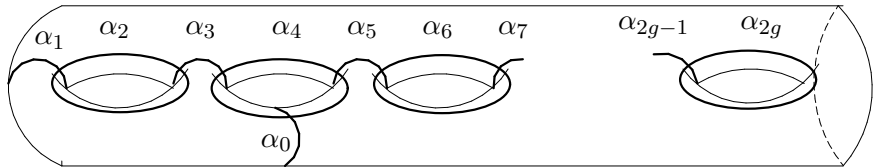
\begin{figure}
\setlength{\unitlength}{0.00056868in}
\begingroup\makeatletter\ifx\SetFigFont\undefined%
\gdef\SetFigFont#1#2#3#4#5{%
  \reset@font\fontsize{#1}{#2pt}%
  \fontfamily{#3}\fontseries{#4}\fontshape{#5}%
  \selectfont}%
\fi\endgroup%
{\renewcommand{\dashlinestretch}{30}
\begin{picture}(8034,1527)(0,-10)
\drawline(420.000,870.000)(493.364,781.874)(579.228,705.875)
	(675.613,643.759)(780.294,596.957)(890.855,566.550)
	(1004.747,553.239)(1119.340,557.332)(1231.992,578.734)
	(1340.103,616.950)(1441.178,671.100)(1532.886,739.934)
	(1613.110,821.864)(1680.000,915.000)
\put(400,1100){\makebox(0,0)[b]{$\alpha_1$}}
\put(1000,1200){\makebox(0,0)[b]{$\alpha_2$}}
\put(1800,1200){\makebox(0,0)[b]{$\alpha_3$}}
\put(2500,150){\makebox(0,0)[b]{$\alpha_0$}}
\put(2600,1200){\makebox(0,0)[b]{$\alpha_4$}}
\put(3300,1200){\makebox(0,0)[b]{$\alpha_5$}}
\put(4000,1200){\makebox(0,0)[b]{$\alpha_6$}}
\put(4700,1200){\makebox(0,0)[b]{$\alpha_7$}}
\put(6000,1200){\makebox(0,0)[b]{$\alpha_{2g-1}$}}
\put(7000,1200){\makebox(0,0)[b]{$\alpha_{2g}$}}\drawline(1545.000,735.000)(1442.369,787.564)(1334.251,827.645)
	(1222.156,854.681)(1107.654,868.294)(992.346,868.294)
	(877.844,854.681)(765.749,827.645)(657.631,787.564)
	(555.000,735.000)
\drawline(3390.000,870.000)(3463.364,781.874)(3549.228,705.875)
	(3645.613,643.759)(3750.294,596.957)(3860.855,566.550)
	(3974.747,553.239)(4089.340,557.332)(4201.992,578.734)
	(4310.103,616.950)(4411.178,671.100)(4502.886,739.934)
	(4583.110,821.864)(4650.000,915.000)
\drawline(4515.000,735.000)(4412.369,787.564)(4304.251,827.645)
	(4192.156,854.681)(4077.654,868.294)(3962.346,868.294)
	(3847.844,854.681)(3735.749,827.645)(3627.631,787.564)
	(3525.000,735.000)
\drawline(6270.000,915.000)(6343.364,826.874)(6429.228,750.875)
	(6525.613,688.759)(6630.294,641.957)(6740.855,611.550)
	(6854.747,598.239)(6969.340,602.332)(7081.992,623.734)
	(7190.103,661.950)(7291.178,716.100)(7382.886,784.934)
	(7463.110,866.864)(7530.000,960.000)
\drawline(7395.000,780.000)(7292.369,832.564)(7184.251,872.645)
	(7072.156,899.681)(6957.654,913.294)(6842.346,913.294)
	(6727.844,899.681)(6615.749,872.645)(6507.631,832.564)
	(6405.000,780.000)
\drawline(1905.000,870.000)(1978.364,781.874)(2064.228,705.875)
	(2160.613,643.759)(2265.294,596.957)(2375.855,566.550)
	(2489.747,553.239)(2604.340,557.332)(2716.992,578.734)
	(2825.103,616.950)(2926.178,671.100)(3017.886,739.934)
	(3098.110,821.864)(3165.000,915.000)
\drawline(3030.000,735.000)(2927.369,787.564)(2819.251,827.645)
	(2707.156,854.681)(2592.654,868.294)(2477.346,868.294)
	(2362.844,854.681)(2250.749,827.645)(2142.631,787.564)
	(2040.000,735.000)
\drawline(510.000,1500.000)(410.503,1450.281)(318.819,1387.308)
	(236.700,1312.287)(165.719,1226.653)(107.232,1132.044)
	(62.359,1030.270)(31.958,923.278)(16.611,813.114)
	(16.611,701.886)(31.958,591.722)(62.359,484.730)
	(107.232,382.956)(165.719,288.347)(236.700,202.713)
	(318.819,127.692)(410.503,64.719)(510.000,15.000)
\drawline(7755.000,15.000)(7827.055,106.784)(7888.890,205.741)
	(7939.797,310.739)(7979.195,420.575)(8006.630,533.993)
	(8021.790,649.692)(8024.501,766.349)(8014.732,882.627)
	(7992.594,997.196)(7958.341,1108.744)(7912.366,1215.993)
	(7855.194,1317.716)(7787.480,1412.748)(7710.000,1500.000)
\dashline{3}(7710.000,1500.000)(7649.004,1403.998)(7597.059,1302.811)
	(7554.603,1197.291)(7521.992,1088.325)(7499.500,976.830)
	(7487.317,863.744)(7485.546,750.017)(7494.200,636.605)
	(7513.208,524.464)(7542.410,414.535)(7581.559,307.744)
	(7630.327,204.989)(7688.303,107.133)(7755.000,15.000)
\thicklines
\drawline(2580.000,15.000)(2662.258,94.863)(2708.384,199.825)
	(2711.584,314.430)(2671.387,421.803)(2593.713,506.131)
	(2490.000,554.999)
\drawline(555.000,735.000)(530.195,846.468)(461.828,937.936)
	(361.945,993.288)(248.144,1002.771)(140.477,964.715)
	(57.914,885.825)(15.001,780.000)
\drawline(3525.000,735.000)(3510.252,844.614)(3449.405,936.975)
	(3354.510,993.788)(3244.362,1003.801)(3140.777,965.032)
	(3064.273,885.158)(3030.000,780.000)
\drawline(4785.000,1005.000)(4675.042,996.581)(4580.973,939.027)
	(4523.419,844.958)(4515.000,735.000)
\drawline(6405.000,780.000)(6356.998,907.083)(6261.887,1004.080)
	(6135.770,1054.568)(6000.000,1050.000)
\drawline(2040.000,735.000)(2025.252,844.614)(1964.405,936.975)
	(1869.510,993.788)(1759.362,1003.801)(1655.777,965.032)
	(1579.273,885.158)(1545.000,780.000)
\put(1050,780){\ellipse{1260}{540}}
\put(2530,733){\ellipse{1260}{540}}
\put(4036,771){\ellipse{1260}{540}}
\put(6875,814){\ellipse{1260}{540}}
\thinlines
\path(510,1500)(7710,1500)(7755,1500)
\path(510,60)(510,15)(7755,15)
\end{picture}
}

\caption{Humphries' Dehn twists generating $\G_{g,1}$}
\end{figure}

Genus three  is important because Wajnryb's presentation of $\G_{g,1}$
shows that all the exotic (non Artin) relations  manifest themselves here. Precisely: if $g\ge 3$ and $\alpha_0,\dots ,\alpha_{2g}$ are the embedded circles in $\Sigma_{g,1}$ as indicated in Figure 2, then following Humphries \cite{hum} the associated Dehn twists $D_{\alpha_0},\dots ,D_{\alpha_{2g}}$ generate $\G_{g,1}$
and a theorem of Wajnryb \cite{wajn} implies that the relations among them follow from the obvious Artin relations recalled above and the exotic relations that involve
$D_{\alpha_0},\dots ,D_{\alpha_6}$ only and define $\G_{3,1}$. This fact is 
to an algebraic geometer rather suggestive: the locus in $\Mcal_{g,\vec{1}}$ for which
the underlying curve is hyperelliptic and the point is a Weierstra\ss\ point is a smooth closed subvariety 
$\Hcal_{g,\vec{1}}$ of $\Mcal_{g,\vec{1}}$. It is an Eilenberg-MacLane space for
the \emph{hyperelliptic mapping class group} of genus $g$: if $\iota$ is any involutional symmetry of
$\Sigma_{g,1}$ with $2g+1$ (interior) fixed points, for example in Figure 2 the symmetry in a horizontal axis, then this is the connected component group of
the group of diffeomorphisms  that leave its boundary pointwise fixed and commute with $\iota$. It is known that this group has a presentation with generators the Dehn twists $D_{\alpha_1},\dots ,D_{\alpha_{2g}}$ subject to the  Artin relations (this is indeed Artin's braid group with $2g+1$ strands). 
Wajnryb's presentation can be interpreted as saying that the inclusion in
$\Mcal_{g,\vec{1}}$  of the union of $\Hcal_{g,\vec{1}}$ and a regular neighborhood of a certain locus in the Deligne-Mumford boundary (that is easily specified) intersected with $\Mcal_{g,\vec{1}}$ induces
an isomorphism on fundamental groups. It is an interesting challenge to 
algebraic geometers  to prove this with methods indigenous to their field (such
as those developed by Zariski and the later refinements by Fulton-Lazarsfeld). At the very least we would then have reproduced the above presentation of $\G_{g,1}$ for 
$g\ge 3$ without the aid of a computer, but most likely this leads also to some conceptual gain.

\section{Orbit spaces of finite reflection groups--a brief review}

In this short section we recall some facts concerning the  Artin groups 
attached to finite reflection groups. Let $W$ be a finite group of
transformations of a real, finite dimensional vector space $V$
which is generated 
by reflections. The fixed point hyperplanes of the reflections in
$W$ decompose $V$ into simplicial cones, called {\it faces}. An open face is
called a {\it chamber}. Let us fix a chamber $C$ and let $\{ H_i\}_{i\in I}$
be the collection of its supporting hyperplanes. Each $H_i$ is the fixed point
set of a reflection $s_i\in W$. If $m_{ij}$ denotes the order of $s_is_j$, 
then let $\G$ be the graph on the index set $I$ that has $m_{ij}-2$ edges
connecting $i$ and $j$ (with $i,j\in I$ distinct). The obvious
homomorphism $W_\G\to W$ is known to be an isomorphism \cite{bourb}
and we will therefore identify $W$ with $W_\G$. The element $w_0\in W$ that
sends $C$ to $-C$ is called the {\it longest element of} $W$. Notice that 
$-w_0$ is a linear transformation of $V$ which maps $C$ to itself and hence
induces a permutation in $I$. Since $w_0^2=1$, this permutation  
is called the {\it canonical involution} of $\G$. 

The Artin group $\Ar_\G$ has now an interpretation as a fundamental group:
The action of $W$ on
\[
\VV^\circ :=\VV-\cup \{ \HH \, :\, H\text{ refl. hyperpl. of } W\}
\]
is free (the use blackboard font indicates that we have complexified the
subspace of $V$ of that same name). Given $x\in V$, let $C_x$ denote
the intersection of halfspaces containing $C$  bounded by a 
reflection hyperplane passing through $x$. So $C_x$ is an open cone
containing $C$ and $C_x=V$ if $x$ is not in any reflection hyperplane.
The subset $\UU\subset \VV$ consisting of the $x+\sqrt{-1}y$ 
with the property
that $y\in C_x$ is open and contained in $\VV^\circ$.  Moreover,
$\UU$ is starlike with respect to any point in $\sqrt{-1}C$ and hence 
contractible. Now fix a base point $*\in C$. Then for every $w\in W$ 
there is a well-defined relative homotopy class of curves in $\UU$ 
of curves connecting $*$ with  $w(*)$. Denote the image of that class 
in the orbit space $\VV^\circ_W$ by $t(w)$ so that 
$t(w)\in\pi_1(\VV^\circ_W, *)$. 

\begin{proposition}[Brieskorn \cite{bries}]\label{finitecoxeter}
The group $\pi_1(\VV^\circ_W, *)$ is generated by the elements $t(s_i)$.
These generators satisfy the Artin relations 
\[ 
t(s_i)t(s_j)t(s_i)\cdots = t(s_j)t(s_i)t(s_j)\cdots 
\text{ ($m_{ij}$ letters on both sides)}
\] 
and the resulting homomorphism  $\Ar_\G\to \pi_1(\VV^\circ_W, *)$ is an 
isomorphism. This homomorphism maps the Garside element $\Delta_\G$ to
$t(w_0)$, where $w_0\in W$ is the element that sends $C$ to $-C$.
\end{proposition}

If we identify $\pi_1(\VV^\circ_W, *)$ with $\Ar_\G$, then the 
map $w\in W_\G\mapsto t(w)\in\Ar_\G$ is in fact
a section $\sigma$ of the natural homomorphism 
$\Ar_\G\to W_\G$. This section is multiplicative: 
$t(ww')=t(w)t(w')$ as long as the length of 
the product $ww'$ equals the sum of the lengths $w$ and $w'$. For instance,
for distinct $i$ and $j$, we have 
$t(s_i)t(s_j)t(s_i)\cdots =t(s_is_js_i\cdots )$ 
($m_{ij}$ letters) and this is also the Garside element of the Artin 
subgroup generated by $t(s_i)$ and $t(s_j)$. 
Such a section exists for any Artin system (wether it is of finite type 
or not).

It is well-known that if $W$ is irreducible in the sense that it 
cannot be written as a product of two nontrivial reflection groups, 
then the center of $W$ is $W\cap \{\pm 1_V\}$. Brieskorn-Saito and Deligne
show that the center of $\Ar_\G$ is then generated by $\Delta_\G$ or 
$\Delta^2_\G$ according to whether or not $-1_V\in W$.

\section{Affine reflection groups}\label{artin:affine}

In this section $A$ is a real, finite dimensional affine space 
with translation space $V$. The latter is a vector
space, but since one usually writes a
transformation group multiplicatively, we write
$\tau_v:A\to A$ for translation over $v\in V$.  We are also given
a group $W$ of affine-linear  transformations of $A$ generated by 
reflections which acts properly discretely with compact orbit space. 
(For a group generated by 
reflections this amounts to its translation subgroup being of finite index
and defining a discrete cocompact lattice in $V$.) 
In order to avoid irrelevant complications, we assume that $W$ acts 
irreducibly in $V$. 

\subsection*{Realization as a Coxeter group}
We may identify $W$ as a Coxeter group in much the same way as for
a finite reflection group. The general reference for this section is
\cite{bourb}. The fixed point hyperplanes of the 
reflections in $W$ decompose $A$ into bounded (relatively open) polyhedra, 
which turn out to be simplices. The 
group $W$ acts simply transitively on the collection of
open simplices (the chambers). We fix one such chamber $C$. Then 
$\bar C$ is a fundamental domain for $W$ in $A$.
Let $\{ H_i\}_i$ be the collection of the supporting hyperplanes of $C$, 
and let  $s_i\in W$ be the reflection that has $H_i$ as its fixed point 
hyperplane. If $m_{ij}$ denotes the order  of $s_is_j$, 
then let $\G$ be the graph on the index set $I$ that 
has $m_{ij}-2$ edges connecting $i$ and $j$ (with $i,j\in I$ distinct). 
Then $\G$ is connected and of affine type and the obvious homomorphism 
$W_\G\to W$ is an isomorphism \cite{bourb}.
We shall therefore identify $W_\G$ with $W$. The Tits construction 
\cite{bourb} shows that any connected graph of affine type so arises.

We denote the translation lattice of $W$ by $Q\subset V$ and the 
image of $W$ in $GL(V)$ by $W^\lin$. Then $W$ is isomorphic to
the semidirect product $Q\rtimes W^\lin$, for there exist points $a\in A$ 
whose stabilizer $W_a$ maps isomorphically onto $W^\lin$. Such a point is 
called the {\it special point}\index{special!point} for this action. 
If $a, b\in A$ are special, then clearly translation over $b-a$ normalizes 
$W$. Conversely, if a translation normalizes $W$, then
it sends any special point  to a special point. So the lattice $P\subset V$
of translations that normalize $W$ acts (simply) transitively on the collection
of special points. This lattice contains $Q$ and is called the 
{\it weight lattice}. (In fact, $W^\lin$ is the Weyl group of a finite, 
irreducible, reduced root system whose span is  $Q$ and $P$ is the 
lattice dual to the span of the dual root system.)
The notions of a special point and of a special vertex are related 
as follows: if for $i\in I$, $a_i$ denotes the unique vertex of $\bar C$ 
that is not contained $H_i$, then $i$ is special precisely when $a_i$ is.

Consider the normalizer of $W$ in the group of affine
transformations of $A$. Then the stabilizer of $C$ in this normalizer
acts faithfully on $C$, and hence likewise on the Coxeter
system $(W, \{ s_i\}_{i\in I})$ and the graph $\G$. The Tits construction
makes evident that the image is
all of $\aut (\G)$ and so the normalizer is in fact 
$W\rtimes \aut (\G)$. The translation 
subgroup of the normalizer is the lattice $P$ and hence 
$\aut (\G)$ contains $Q/P$ as a distinguished subgroup. Now
the Corollary to Prop.\ 6 of Ch.~VI, \S 2 of \cite{bourb} can be 
stated as follows:

\begin{lemma} 
The group $P/Q$ acts simply transitively on $I_\spec$. 
\end{lemma}

The characterization of $\aut(\G )$ as the stabilizer of $C$ in
the normalizer of $W$ leads to a geometric interpretation of the 
elements $g_i$ and $g_{ij}$ of $\G$ (with $i$ and $j$ special) 
defined earlier. For every (nonempty) face $F$ of $C$ we 
denote by $W_F$ the subgroup of $W$ that leave $F$ pointwise fixed.
This group is generated by the $s_i$, $i\in I$, with 
$H_i\supset F$ and $(W_F,\{ s_i\}_{H_i\supset F})$ is a Coxeter system
of finite type. Let $w_F$ denote the element of
$W_F$ that maps $C$ to the chamber that is opposite $F$.
This is an element of order two and is the `longest'
element of $W_F$. In case $F$ consists of a single vertex $a_i$ or
is the span of two such vertices $a_i,a_j$, then we also write $w_i$ 
resp.\ $w_{ij}$ instead of $W_F$. 
Denote by $\iota_j$ resp.\ $\iota_{ij}$ denote the
central symmetry of $A$ in the vertex $a_j$ resp.\
the midpoint of $a_i$ and $a_j$: $\iota_j(a) =a-2(a-a_j)$ and
$\iota_{ij}(a) =a-(2a-a_i-a_j)$. 

\begin{lemma}\label{vertexinvolution} 
Assume that $j\in I$ has the property that
$\iota_j$ normalizes $W$. Then $w_j\iota_j$
is an involution $g_j$ of $C$ whose image in $\aut(\G )$ fixes
$j$. It preserves each connected component of $\G_j$ and acts there 
as the canonical involution. The assumption is
satisfied when $j$ is special (and in that case $\G_j$ is connected, 
so that $g_j$ acts on $\G_j$ as the canonical involution). 
\end{lemma} 
\begin{proof} The first two statements are clear from the definitions.
If $j\in I$ is special, then $W=Q\rtimes W_j$ and so $\iota_j$
normalizes $W$. 
\end{proof}

It may happen that  $\iota_j$ normalizes $W$ for
some nonspecial $j\in I$. An example that will appear here
is when $\G$ is of type $\hat E_7$ and $j$ is the end of
the short branch.

If two special points of $A$ belong to the same chamber
(i.e., if the open interval they span is a one dimensional
face), then in accordance with received terminology, 
their difference in $V$ is called a {\it
minuscule weight}\index{minuscule!weight}, reason for us
to refer to the corresponding pair of special points as 
a {\it minuscule pair}\index{minuscule!pair}.
Since $\iota_{ij}\iota_i$ is translation over $a_j-a_i$, 
the involution $\iota_{ij}$ normalizes $W$.
Define $g_{ij}$ by the property that $\iota_{ij}=w_{ij}g_{ij}$.
This is an affine linear transformation of $A$ which preserves $C$.
 
Recall that $S(\G )$ stands for the subgroup 
of $\aut(\G)$ generated by the
elements $g_i$  and $g_{ij}$ with $i$ special resp.\ $i,j$
a minuscule pair. We shall write 
\[
\sw:=W_\G \rtimes S(\G ) \text{ and } \sar:=\Ar_\G \rtimes S(\G ). 
\]
It is clear that $\sw$ can be obtained from $\sar$ 
by imposing the relations $t_i^2\equiv 1$ (all $i$). We put
\[
\swbar:=\{ \pm 1_V\}.W^\lin \subset\GL (V), \quad \wbar:=
\{ \pm 1_V\}.W^\lin / \{ \pm 1_V\}.
\]
So the first group or the second group equals $W^\lin$, depending on
whether or not $-1_V\in W^\lin$. We justify this notation  
by the following result.

\begin{corollary}\label{cor:wreduction}
The group $S(\G)$ contains $P/Q$ as a normal subgroup. 
The quotient is cyclic of order two
if $W^\lin$ does not contain $-1_V$ and is trivial otherwise,
in other words, the natural homomorphism $\sw\to \GL (V)$ has
kernel $P$ and image $\swbar$.
Moreover, $\wbar$ is obtained as 
a quotient of the reduced  Artin group $\arbar_\G$ by imposing the 
relations $t_i^2\equiv 1$ (all $i$).
\end{corollary}
\begin{proof} We use the fact $\sw$ is generated by $W$ and the central 
symmetries $\iota_i$ and $\iota_{ij}$. 
We have seen that $\sw$ contains the translation over $a_j-a_i$ if
$(i,j)$ is a minuscule pair. So $\sw$ contains $P$ as a translation 
subgroup. Since  $P$ is the full group of translations normalizing $W$,
it follows that $P$ is the kernel of the natural homomorphism 
$\sw\to \GL (V)$. Since $\iota_i$ and $\iota_{ij}$ act in $V$ 
as minus the identity it follows that the image is 
$\sw^\lin$, as asserted.

Finally, the set of relations defining $\arbar$ yield
in $\sw$ the relations  $w_ig_i\equiv 1$ ($i$ special) and 
$w_ig_i \equiv w_{ij}g_{ij}$ ($(i,j)$ a minuscule pair). 
Since $w_ig_i=\iota_i$ and $w_{ij}g_{ij}=\iota_{ij}$,
we see that imposing these relations amount to killing $P$ 
(so that are dealing with quotient of $\swbar$) 
and subsequently the common image $-1_V$ of the central symmetries.
The result is indeed $\wbar$.
\end{proof}

Since conjugation by $g_{ij}$ interchanges $g_i$ and $g_j$, the $g_i$'s
make up a single conjugacy class in $S(\G )$. So if $P/Q$
is cyclic and $W^\lin$ does not contain $-1_V$, then $S(\G)$ 
has the structure of a dihedral group whose group of rotations 
is $P/Q$ and whose reflections are the $g_i$'s and the $g_{ij}$'s.

Let us go through the list of connected graphs of affine type and see
what we get (see for instance \cite{bourb} for a classification).
 When there is just one special point ($\hat
E_8,\hat F_4,\hat G_2$) or two ($\hat A_1$, $\hat B_{l\ge 2}$, 
$\hat C_{l\ge 3}$, $\hat E_7$), there is
little to say: in these cases, the involutions $g_i$ must be 
trivial and $S(\G)=\aut (\G)$ (and of order at most two).

In case $\hat D_{l\ge 4}$ the special
vertices are the four end vertices. When $l\ge 5$
the graph $\hat D_l$ has two branch points, each connected
with two end vertices. When $l$ is odd, then we arrange the
four end vertices as the vertices of a square such that
vertices separated by a single branch vertex are opposite and
$S(\G )$ can be identified with the corresponding dihedral
group of order 8. When $l$ is even, the involutions $g_i$
are all trivial and so $S(\G )=P/Q$. The latter is
isomorphic to $\ZZ/2\oplus \ZZ/2$.

The remaining cases are $\hat A_{l\ge 2}$ and $\hat E_6$.
In the first case, $\G$ is an $(l+1)$-gon all of whose
points are special and $S(\G )=\aut (\G)$ is a dihedral
group of order $2(l+1)$.
In case $\hat E_6$ the situation is basically the same: we
have three special points, the $g_i$'s are the
transpositions and $g_i=g_{jk}$ if $i,j,k$ are 
special and mutually distinct. 

From this description we deduce that  $S(\G )=\aut (\G )$
and is a dihedral group, unless $\G$ is of type $\hat
D_\even$, in which case $S(\G )\cong \ZZ/2\oplus \ZZ/2$.

\subsection*{Affine Artin groups as fundamental groups}
We continue with the situation of previous section. The complexification
$\AA$ of the affine space $A$ is naturally isomorphic to $A\times V$ 
(the affine transformation group of $A$ acts after complexification 
on the imaginary part through $GL(V)$). 

Let $\AA^\circ\subset\AA$ be the complement of the union of all the 
complexified reflection hyperplanes of $W$. Then $W$ acts freely
on $\AA^\circ$ so that its orbit space $\AA^\circ_W$ is a
complex manifold. Its fundamental group can be determined in much the same
way as in the finite reflection group case. We begin with defining
a contractible $W$-open neighborhood $\UU$ of the $W$-orbit of 
$C\times\{ 0\}$ in $\AA^\circ$. Given $x\in A$, then
the intersection of halfspaces containing $C$  bounded by a 
reflection hyperplane passing through $x$ is of the form 
$x+C_x$ for some open 
linear convex cone $C_x$ in $V$. We now take $\UU$ to be the set of 
$(x,y)\in\AA$ with $y\in C_x$. This is indeed an open subset of $\AA^\circ$. 
The projection $\UU\to A$ is surjective and the preimage of the
star of every face of $A$ is starlike (see the argument used
for a finite reflection group) and hence contractible.
This implies that $\UU$ is contractible. 

Let $x_o\in C$ be the barycenter (and hence fixed under 
any automorphism of $C$) and let $*:=(x_o,0)\in \UU$. For every
$w\in W$, there is unique homotopy class of curves in $\UU$ from $*$ to
$w(*)$; we denote the image of that homotopy class in $\pi_1(\AA^\circ_W,*)$
by $t(w)$. In analogy to \ref{finitecoxeter} we have:

\begin{proposition}[Nguy\^e\~ n Vi\^et D\~ ung \cite{nguy}] \label{affinecoxeter}
The group $\pi_1(\AA^\circ_W, *)$ is generated by the elements $t(s_i)$.
These generators satisfy the Artin relations 
\[ 
t(s_i)t(s_j)t(s_i)\cdots = t(s_j)t(s_i)t(s_j)\cdots 
\text{ ($m_{ij}$ letters on both sides)}
\] 
and the resulting homomorphism  $\Ar_\G\to \pi_1(\AA^\circ_W, *)$ is an 
isomorphism. If we identify these two groups by means of this isomorphism, 
then $w\in W_\G\mapsto t(w)\in\Ar_\G$ is a section of the projection 
$\Ar_\G\to W_\G$ with the property that $t(ww')=t(w)t(w')$ if
the length of $ww'$ is the sum of the length of $w$ and of $w'$.
\end{proposition}

The group $\sw$ acts also on $\AA^\circ$ and we thus find

\begin{corollary}\label{orbifold:semidirect}
The orbifold fundamental group of $\AA^\circ_\sw$ is $\sar$.
\end{corollary}

Under this identification (as defined in the introduction of this paper), 
an element $g\in S(\G )$ is represented by the pair consisting of the 
constant path $[*]$ at $*$ 
and the element $g$. So $t(w).g$ is represented by the pair 
$(t(w),wg)$. On the other hand, $gt(w)$ is represented by $(t({}^gw), gw)$,
We extend the section $W\to \Ar$ to a section 
$t$ of $\sar\to \sw$ by putting $t(wg):= t(w)g$. 

\begin{lemma}\label{minuscule:secidentity}
If $(i,j)$ is a minuscule pair, then 
$t(\tau_{a_j-a_i})=(\Delta_{ij}g_{ij})^{-1}\Delta_ig_i$.  
\end{lemma}
\begin{proof} From 
$\tau_{a_i-a_j}=\iota_{ij}\iota_i=w_{ij}g_{ij}w_ig_i
= w_{ij}w_jg_{ij}g_i=w_{ij}^{-1}w_jg_{ij}g_i$ we see that
$t(\tau_{a_i-a_j})=t(w_{ij}^{-1}w_j)g_{ij}g_i$.
Since the length of $w_j$ is the sum of the lengths of 
$w_{ij}$ and $w_{ij}^{-1}w_j$, we have 
$\Delta_j=\Delta_{ij}t(w_{ij}^{-1}w_j)$. Substituting this
in the preceding gives
$t(\tau_{a_i-a_j})=\Delta_{ij}^{-1}\Delta_jg_{ij}g_i=
\Delta_{ij}^{-1}g_{ij}\Delta_ig_i=(\Delta_{ij}g_{ij})^{-1}\Delta_ig_i$. 
\end{proof} 

\subsection*{Toric structure of the orbifold}
We can also arrive at $\AA^\circ_\sw$ by starting from the 
intermediate orbit space $\AA_P$.  Notice that  
$\AA_P$ is naturally an algebraic torus whose identity element is the
image of the special orbit. It is clear $\AA_P$ is as a torus 
canonically isomorphic the `adjoint' torus $\CC^\times\otimes_\ZZ P$,
where the adjective canonically guarantees (via Corollary \ref{cor:wreduction})
equivariance relative to $\swbar=SW /P=\{\pm 1_V\}\cdot W^\lin$.
We shall write $\TT$ for $\CC^\times\otimes_\ZZ P$.  
The open subset $\TT^\circ$ 
is the complement of the union of the reflection hypertori 
(relative to the $\swbar$-action). We may also consider the
orbit space $\TT_{\pm 1}$ and its open subset $\TT^\circ_{\pm 1}$.
Note that $\wbar$ acts faithfully on both. It is clear that
\[
\AA^\circ_{\sw}\cong \TT^\circ_{\swbar}.
\]
We shall use this identification to describe 
certain extensions of $\AA^\circ_{\sw}$ as an orbifold by means of 
$\swbar$-equivariant extensions of $\TT^\circ$.

We take the occasion to prove Proposition \ref{prop:smallgroup}.

\begin{proof}[Proof of Proposition \ref{prop:smallgroup}]
In view of Corollary \ref{cor:wreduction} it suffices to prove that the square
of every generator of $\hat A_l$ becomes 1 in $\overline{W}_{\hat A_l}$. 

To this end we label the vertices of $\hat A_l$ in cyclic order by $\ZZ/(l+1)$.
If $\zeta$ denotes the rotation of the graph that sends the vertex with index $i$ to the one with index $i+1$, then the relations  imposed on $\overline{Ar}_{\hat A_l}$
imply that
$\Delta_{\ZZ/(l+1)-\{ i\}}\equiv\Delta_{\ZZ/(l+1) -\{ i,j\}}\zeta^{i-j}$
for all distinct pairs $(i,j)$ in $\ZZ /(l+1)$.
If we keep in mind that
\[
\Delta_{A_k}=(t_1t_2\cdots t_k)(t_1t_2\cdots t_{k-1})\cdots (t_1t_2)t_1,
\]
then we see that for $(i,j)=(0,l)$ this gives $t_1t_2\cdots t_l\equiv
\zeta$ and
for $(i,j)=(0,l-1)$  $t_1t_2\cdots t_l t_1t_2\cdots t_{l-1}\equiv
t_l\zeta^2$.
Together these identities imply $t_l^2\equiv 1$. Any other generator
$t_i$ is conjugate to $t_l$ and so we also have that $t_i^2\equiv 1$.
\end{proof}

\subsection*{Effect of a divisorial extension on the fundamental group}
Let $X$ be a connected complex manifold, 
and $Y\subset X$ a smooth hypersurface and  write $U$ for
its complement $X-Y$. 
We fix a base point $*\in U$. Then the natural homomorphism
$\pi_1(U,*)\to \pi_1(X,*)$ is surjective with with its kernel
normally generated by the simple (positively  
oriented) loops around $Y$.  Two such loops are conjugate in $\pi_1(U,*)$
if they encircle the same connected component. So we may say that
$\pi_1(X,*)$ is obtained from $\pi_1(U,*)$ by imposing 
a relation $\rho_i\equiv 1$ for every connected component 
$Y_i$ of $Y$ (with $\rho_i$ represented by a simple loop around $Y_i$).

This is still true in an orbifold setting: if a 
group $G$ acts properly discontinuously on $Y$ and preserves $H$, 
then the natural homomorphism
\[
\pi_1^\orb(U_G,*)\to \pi_1^\orb(X_G,*)
\]
is surjective with its kernel normally generated by the $\rho_i$'s, 
where $i$ now runs over a system of representatives
of the $G$-action on $\pi_0(H)$.

\subsection*{Toric extensions}
Any indivisible $p\in P$ defines an injective homomorphism
$\rho_p :\CC^\times\to \TT$. This defines a smooth
partial torus embedding $\TT\subset \TT\times^{\rho_p}\CC$
which adds to $\TT$ a divisor $\TT (p)$ canonically isomorphic 
to the cokernel of $\rho_p$. This can be done independently (and
simultaneously) for any finite set $\Sigma\subset P$ of indivisible
vectors. The corresponding partial torus embedding, 
denoted  $\TT_\Sigma$, is smooth. If 
$\Sigma$ is $\swbar$-invariant, then we can form the 
orbifold $(\TT_\Sigma)_{\swbar}$.

\begin{lemma}\label{toric} 
Let $p\in P$ be indivisible and nonzero. 
Then the inclusion of orbifolds
$\TT^\circ_\swbar\subset (\TT_{\swbar .p}^\circ)_\swbar$ 
introduces on orbifold fundamental groups the relation 
$t(\tau_{-p})\equiv 1$ in $\sar$.
\end{lemma} 
\begin{proof} 
Let $x_o\in C$ be as chosen earlier and let $y\in V$ be not in any
reflection hyperplane. Consider 
the path $\omega_y$ in $\AA$ from $*=(x_o,0)$ to $\tau_{p}(*)=(x_o-p,0)$ that 
traverses the segments successively connecting 
$(x_o,0)$, $(x_o,y)$, $(x_o-p,y)$ and $(x_o-p,0)$. 
It need not lie in $\AA^\circ$, but if we
replace $y$ by $y+\lambda p$ with $\lambda$ large, it will even lie in
$\UU$ (as defined in \ref{affinecoxeter}) and will hence represent $t(\tau_{-p})$. 
The image of this path in $\TT$ defines a simple negatively oriented loop
around $\TT(p)$. 
\end{proof}

An extension of $\TT^\circ_\swbar$ as above will be refered
to as a \emph{toric extension}.\index{extension!toric}
Lemmas \ref{toric} and \ref{minuscule:secidentity} imply:

\begin{corollary}\label{toriccor}
The toric extension of $\TT^\circ_\swbar$ defined by the set of
minuscule weights has orbifold fundamental group canonically isomorphic
to the quotient of the extended Artin group $\Ar\rtimes S(\G)$ by the relations 
$\Delta_{ij}g_{ij}=\Delta_ig_i$ with $i,j$ special and distinct.
\end{corollary}

\subsection*{Blowup extensions}
We next consider certain equivariant blowups. 
The image of $A$ in $\TT$ is a compact 
torus and the decomposition of $A$ in faces gives a similar decomposition of
that torus. This decomposition is finite; in particular, 
we have finitely many vertices (of finite order).  
A $\swbar$-invariant set of vertices in $\TT$ determines an
$S(\G )$-invariant subset of $I$ and vice versa. 

Fix a $S(\G )$-invariant subset $J\subset I$, let $v(J)$ denote 
the corresponding $\swbar$-invariant set of vertices 
in $\TT$ and consider the blowup of $\TT$ at $v(J)$: 
\[
\bl_{v(J)}(\TT)\to \TT .
\]
The exceptional set $E_j$ over the image of
$a_j$, $j\in J$, in $\TT$ can be
identified with the projective space $\PP (\VV)$. This
exceptional set is of course invariant under $\sw_{a_j}$ and the
identification with $\PP (\VV)$ is equivariant. So the
intersection of $\bl_{v(J)}(\TT)^\circ$ with $E_j$ 
can be identified with the complement of the
union of reflection hyperplanes of $W_{\G_j}$ in $\PP
(V)$. Since the finite group $\swbar$ acts on $\bl_{v(J)}(\TT)^\circ$
we can form $\bl_{v(J)}(\TT)^\circ_\swbar$; it contains  $\TT^\circ_\swbar$
as the complement of a closed hypersurface.

\begin{lemma}\label{blowup} 
Suppose that $j\in I$ is special, or more generally, 
suppose that $j\in I$ is such that $\iota_j\in\sw$, 
so that we can define $g_j\in S(\G)$ by the proprerty $w_jg_j=\iota_j$.
Then the open immersion of orbifolds
$\TT^\circ_\swbar\subset \bl_{v(J)}(\TT)^\circ_\swbar$ 
introduces on orbifold fundamental groups the relation 
$\Delta_{\G_j}g_j\equiv 1$ in $\sar$.
\end{lemma} 
\begin{proof} An element of $\swbar$
which leaves $E_j$ pointwise fixed
must lie in $\sw_{a_j}$ and induce a
homothety in $V$. Such element must be the central symmetry
relative to $v(j)$ (or the identity). 
So the relation in the orbifold
fundamental group defined by $E_j$ comes from a half loop 
in $\TT^\circ$ around $v(j)$ in $\TT$. This loop is
represented in $\Ar$ by $\Delta_{\G_j}g_j$.
\end{proof}

It is clear that toric and blowup extensions can be performed independently. 
Our main interest concerns the case 
when $J=I_\spec$ (so that $v(J)=\{ 1\})$)
and $\Sigma$ is the set of all minuscule weights. 
In that case we abbreviate
\[
\hat\TT:= \bl_{\{ 1\}}(\TT)(\Sigma),
\]
so that $\hat\TT^\circ \supset \TT^\circ$. This yields an inclusion of orbifolds
$\hat\TT^\circ_\swbar$ and  
from  \ref{toriccor} and \ref{blowup} we immediately get:

\begin{corollary}\label{cor:caseofinterest} 
The  open embedding of orbifolds $\TT_\swbar^\circ\subset \hat\TT^\circ_\swbar$ 
yields on orbifold fundamental groups the
epimorphism $\sar_\G\to\arbar_\G$.
\end{corollary}

In case $\G$ is of type $\hat E_7$, 
we will also consider the case when
$J$ consists of the three end vertices, denoted here
$s,s',n$, where $s,s'$ are the special vertices. Then $n$ satisfies
the hypotheses of Lemma \ref{blowup} with $g_n=g_{s,s'}$ and
we thus obtain the group appearing in the statement 
of Theorem \ref{thm:genus3} as orbifold fundamental group:

\begin{corollary}\label{cor:caseofinterest2} 
Let $\tilde\TT$ be obtained from $\TT$ 
by taking for the toric data the set of all minuscule weights
and for the the blownup data the set of the three end vertices.
Then the open embedding of orbifolds $\TT_\swbar^\circ\subset\tilde\TT^\circ_\swbar$ 
yields on orbifold fundamental groups the introduction of the relation
$\Delta_{A_7}\equiv\Delta_{E_6}$ in $\sar_\G$.
\end{corollary}

\section{Moduli spaces of marked Del Pezzo surfaces}\label{sect:dpmoduli}

Let $S$ be a complex Del Pezzo surface. By definition this
means that $S$ is a complete smooth surface whose
anticanonical bundle $\omega_S^{-1}$ is ample. 
\index{Del Pezzo surface} Its {\it
degree} $d=\omega_S\cdot\omega_S$ is then a positive integer
$\le 8$. The complete linear system $|\omega_S^{-1}|$ is
$d$-dimensional, and its general member
is a genus one curve. It has no base points when $d\ge 2$,
and when  $d\ge 3$ it defines an embedding. 
For $d=3$ the image is of this map is a cubic
surface and for $d=4$ a complete intersection of two
quadrics.  When $d=2$, the map $|\omega_S^{-1}|$ is a
twofold covering of a projective plane which ramifies along
a smooth curve of degree four and this defines a natural 
involution of $S$. Conversely, any smooth
cubic surface, complete intersection of two quadrics or
double cover of a projective plane ramifying along smooth
quartic curve is a Del Pezzo surface (of degree $4,3,2$
respectively). A Del Pezzo surface of degree $d$ is
isomorphic to the blowup of $r:=9-d$ points in a projective
plane in general position  or (when $d=8$) to a product two
projective lines. So there is just one isomorphism type for
$d=5,6,7,9$ and two for $d=8$. The classes $e$ of the
exceptional curves of the first kind are characterized by
the properties $e\cdot e=e\cdot \omega_S=-1$. 
We recall the structure of the Picard group of a Del Pezzo surface.

\subsection*{Del Pezzo lattices}
Let $r$ be an integer $r\ge 0$ and denote the standard
basis vectors of $\ZZ^{r+1}$ by $l,e_1,\dots ,e_r$. We
define an inner product on  $\ZZ^{r+1}$ by requiring that
the basis is orthogonal and $l\cdot l=1$ and $e_i\cdot e_i=-1$. This
is the standard Lobatchevski lattice of rank $r+1$; we
denote it by $\Lambda_{1,r}$.

Following Manin \cite{manin}, we put $k:= 3l-e_1-\cdots
-e_r\in\Lambda_{1,r}$. So $k\cdot k= 9-r$. The orthogonal
complement $Q_r$ of $k$ in $\Lambda_{1,r}$ is negative
definite precisely when $r\le 8$. Then the elements
$\alpha\in Q_r$ with $\alpha\cdot\alpha =-2$ define a (finite)
root system $R_r$; the reflection $s_\alpha$ associated to
$\alpha\in R_r$ is of course orthogonal reflection. We
denote its Weyl group by $W_r$. 

The set of \emph{exceptional
vectors}, $\Ecal_r$, is the  set of $e\in\Lambda_{1,r}$ with
$-e\cdot e=e\cdot k=1$. This set is $W_r$-invariant and consists
of the elements $e_i$, $l-e_i-e_j$, $2l-\sum_{i\in I} e_i$ 
with $I$ a $5$-element subset of $\{ 1,\dots ,r\}$ (so this occurs
only when $r\ge 5$) and
$3l-2e_j-\sum_{i\in I} e_i$ with $I$ a $6$-element subset of 
$\{ 1,\dots ,r\}-\{ j\}$ (occuring only when $r\ge 7$).  

For $r\ge 3$, a root basis of $R_r$ is
$\alpha_1:=l-e_1-e_2-e_3,
\alpha_2:=e_1-e_2,\dots\alpha_r:=e_{r-1}-e_r$. This is also
a basis for $Q_r$ and so $Q_r$  may then be regarded as the
root lattice of $R_r$. The root system $R_r$ is for
$r=3,\dots ,8$ of type $A_1\sqcup A_2$, $A_4$, $D_5$,
$E_6$, $E_7$, $E_8$ respectively. (For $r=2$,
$R_r=\{\pm\alpha_2\}$ and for $r\le 1$, $R_r$ is empty.)
Notice that the subgroup of $W_r$ generated by the roots 
$\alpha_2=e_1-e_2,\dots\alpha_r=e_{r-1}-e_r$ is of type $A_{r-1}$
and can be identified with the symmetric group of the $e_i$'s. It is
also the $W_r$-stabilizer of $l$. We shall concentrate on the cases
$4\le r\le 8$.

\subsection*{Markings} 
Let $S$ be a Del Pezzo surface of degree $d\le 5$ and put $r:=9-d$.
It is clear that there exists an isometry 
$f: \Lambda_{1,r}\to\pic (S)$ which maps $k$ onto
$\omega_S^{-1}$. Such an isometry will be called a \emph{marking}
of $S$. 
Since two such isometries differ by an element of
$W_r$, the latter group permutes simply transitively the
markings of $S$. Notice that a marking $f$ maps an
exceptional vector to the class of exceptional curve of the
first kind and vice versa. It is also known (and again,
easy to verify) that $\alpha\in Q_r$ is a root (i.e.,
satisfies $\alpha\cdot\alpha=-2$) if and only if it is the 
difference of
two mutually perpendicular exceptional vectors, or
equivalently, if and only  if $f(\alpha )$ can be represented by
$E-E'$, where $E$ and $E'$ are disjoint exceptional curves of
the first kind.

The complete linear system $|f(l)|$ is two dimensional
without base points and defines a birational morphism to a
projective plane. Its contracts $r$ exceptional (disjoint)
curves whose classes are $f(e_1),\dots ,f(e_r)$. Their
images will be $r$ points in general position in the sense
that no three are on a line, no six are on a conic and no
eight lie on a cubic which is singular at one of these
points. The linear system $|\omega_S^{-1}|$ is mapped onto
the linear system of cubics passing through these points.

Conversely, if we blow up $r$ points in in
a projective plane in general position we get a
(canonically) marked Del Pezzo surface of degree $d$.
If $d\le 5$, then we can take the first $4$ of these $r=9-d$ points to be 
$[1:0:0]$, $[0:1:0]$, $[0:0:1]$, $[1:1:1]$, so that an open subset of
$(\PP^2)^{5-d}$ is (the base of) a fine moduli space
of marked Del Pezzo surfaces of degree $d$. The group $W_{9-d}$ acts on this open
set and the associated orbifold is a coarse moduli space for 
Del Pezzo surfaces of that degree.

\begin{definition}
A \emph{Del Pezzo pair of degree} $d$ consists of a 
Del Pezzo surface $S$ of degree $d$ and a \emph{singular} anticanonical curve
$K$ of $S$. If also is given a singular point $p$ of $K$, 
then we call $(S,K,p)$ a \emph{Del Pezzo triple of degree} $d$.
\end{definition}

It is  clear that a Del Pezzo pair can be extended to a Del Pezzo triple and
that this can be done in only a finite number of ways.

If $S$ is a Del Pezzo surface of degree $d$ and $p\in S$,
then denote by $L_p$ the linear subsystem of
$|\omega_S^{-1}|$ of anticanonical curves which have a
singular point at $p$. For $d\ge 3$, this subsystem has
dimension $d-3$ always (it is the projective space of hyperplanes
in the receiving projective space of its canonical embedding that
contain the tangent plane of  $S$ at $p$). For $d=2$, $L_p$ is 
nonempty if and only if $p$ lies on the fixed curve (and is then
a singleton and defined by the tangent line of the fixed curve at $p$)'
For $d=1$, generically 12 points $p$ of $S$ have nonempty $L_p$.  

A {\it marking} of a Del Pezzo pair resp.\ triple will be simply a marking of
the first item. There is an evident fine moduli space  $\tilde M(d)$ of  
marked Del Pezzo triples of degree $d$. It is easy to see that  $\tilde M(d)$
is smooth when $d\le 2$. We shall also find that this is so when $d=1$.
The evident action of the group $W_r$ on $\tilde M(d)$ therefore defines
an orbifold $M(d)$, that we interpret as  the coarse
moduli space of such triples. We can similarly define a moduli space for
Del Pezzo pairs (though we will not introduce a notation for it) and 
it is easy to see that the forgetful map from $M(d)$ to this moduli space
is a normalization. Let us record these remarks as a

\begin{lemma} 
For $d=3,4,5$, $M(d)$ is a $\PP^{d-3}$ bundle
over the universal Del Pezzo surface of degree $d$, $M(2)$
is the universal smooth nonhyperelliptic curve of genus
$3$ and $M(1)$ is a 12-fold (ramified) covering of the
coarse moduli space of Del Pezzo surfaces of degree 1.
\end{lemma}

We make the connection  with the setting of Section
\ref{artin:affine} although
the affine Coxeter group will not appear here in a natural manner.
We introduce for $4\le r\le 8$ the affine transformation group 
$\hat{W}_r$ generated by $W_r$ and the translations 
over $Q_r$. This group satisfies the assumptions 
at the beginning of Section \ref{artin:affine}. The quadratic form
on $Q_r$ identifies $Q_r$ with a sublattice of 
$P_r:=\Hom(Q_r,\ZZ)$. The latter is then the group of 
special points of $\Hom(Q_r,\RR)$ and
\[
\TT_r:=\Hom (Q_r,\CC^\times)=\CC^\times\otimes P_r
\]
is the `adjoint torus'. We extend the 
notation that was introduced there in a logical way to the present
situation: we put 
\[
SW_r:=\{\pm 1\}. W_r \text{ and } \wbar_r:=\{\pm 1\}. W_r/\{\pm 1\}.
\] 
Now we are in the setting of 
Section \ref{artin:affine} with $W^\lin =W_r$. 
In particular we have defined $\hat\TT^\circ_r$ and in case $r=7$ 
($W_7$ is of type $E_7$) we also have $(\tilde\TT_7)^{\circ}$.
The following theorem makes the connection between the results of
Section \ref{artin:affine} and some of the theorems stated in Section \ref{sect:results}.

\begin{theorem}\label{thm:moduli}
For $d\le 5$ there is a $\wbar_r$-equivariant open embedding of 
$(\hat\TT^\circ_r)_{\{ \pm 1\}}$ (recall that $r=9-d$) in the moduli space 
$\tilde M(d)$ of marked Del Pezzo 
triples of degree $d$ whose complement is of codimension $\ge 2$ so that
the corresponding open embedding of $(\hat\TT^\circ_r)_{SW_r}$ in $M(9-r)$ induces
an isomorphism between the orbifold fundamental group of
$(\hat\TT^\circ_r)_{SW_r}$ and the reduced orbifold fundamental group of $M(d)$.
The strata of $(\hat\TT^\circ_r)_{SW_r}$ record the Kodaira type of the 
anticanonical curve: the open stratum $(\TT_r)^\circ_{SW_r}$ parameterizes the 
Del Pezzo triples $(S,K,p)$ for which $K$ is an irreducible rational
curve with a simple node at $p$ (Kodaira type $\Ione$), the blowup stratum  
those for which $K$ is a cuspidal curve  
(Kodaira type II) and the toric stratum (which is empty for $r=8$) 
those for which $K$ is  of Kodaira type $\Itwo$. 
\end{theorem}

From this we can deduce two of our main theorems. 

\begin{proof}[Proof that Theorem \ref{thm:moduli} implies  Theorems \ref{thm:main}
and \ref{thm:main2}] This follows  from Theorem \ref{thm:moduli} and  Corollary \ref{cor:caseofinterest}. 
\end{proof}

In case $r=7$, $M(2)$ can be identified with the coarse moduli space
of smooth pointed nonhyperelliptic genus three curves.  

\begin{theorem}\label{thm:moduligenus3}
The open embedding of orbifolds $(\hat\TT_7)^\circ_{SW_7}$ in $M(2)$ extends
to an open embedding of the orbifolds $(\tilde \TT_7)^{\circ}_{SW_7}$ in the 
coarse moduli space $\Mcal_{3,1}$ of pointed smooth genus three curves with
the new stratum (of projective type $A_7$) parameterizing the hyperelliptic curves
pointed with a non-Weierstra\ss\ point. The complement of this immersion
is everywhere of codimension $\ge 2$ and so the immersion induces
an isomorphism of the orbifold fundamental group of $(\tilde\TT_7)^\circ_{SW_r}$
onto the pointed mapping class group $\G_{3,1}$.
\end{theorem}

Another main theorem follows from this:

\begin{proof}[Proof that Theorem \ref{thm:moduligenus3} implies  Theorem \ref{thm:genus3}] Just combine Theorem \ref{thm:moduligenus3} with  Corollary \ref{cor:caseofinterest2}. 
\end{proof}

\subsection*{Reduced Kodaira curves}We use the following defining property.

\begin{definition}
A \emph{reduced Kodaira curve} is a complete connected curve 
of arithmetic genus one whose automorphism group acts 
transitively on its regular part. 
\end{definition}

Let us go through some simple properties of such a curve $K$. 
First of all, if $K$ is singular, then 
\begin{enumerate}
\item[\it{(add)}] $K$ is simply connected and has exactly
one singular point and vanishing first Betti number; it is
then of Kodaira type II (a rational curve with a cusp),
III (two smooth rational curves meeting with multiplicity two)
or IV (isomorphic to the union of three concurrent lines in the plane), or 

\item[\it{(mult)}] $K$ has first Betti number one; it is
then of Kodaira type $I_{n\ge 1}$ (for $n\ge 2$ this is an $n$-gon  consisting of 
smooth rational curves; for $n=1$ this means a rational 
curve with a node). 
\end{enumerate} We then say that $K$ is of {\it additive}
resp.\ {\it multiplicative} type.\index{additive type}
\index{multiplicative type}

For $K$ as above we have the standard short exact sequence
\[
1\to \pic(K)^0\to \pic (K)\to H^2(K;\ZZ)\to 0,
\]
and $\pic(K)^0$ is isomorphic to the additive group $\GG_a$ or
the multiplicative group $\GG_m$ according to its type. 
More precisely, the abelian group $H^2(K;\ZZ)$ is freely generated by the 
irreducible components of $K$ and if $C$ is one such component, then
assigning to $p\in K_\reg\cap C$ the class of $(p)$ defines 
an isomorphism of $K_\reg\cap C$ on the preimage of $[C]$ in $\pic (K)$.  
This shows at the same time that each connected component of $K_\reg$ 
is a principal homogeneous space (torsor) of $\pic(K)^0$. 
Thus $\pic(K)^0$ appears as a  normal subgroup of $\aut (K)$. 
We denote by $G_K$ the quotient group.  
 
If $K$ is of type II, then $K_\reg$ is an affine
line, and $\pic(K)^0$ appears as its group of translations.
So $G_K$ is then isomorphic to $\GG_m$. If $K$ is of
type $\Ione$, then $G_K$ is of order two. In fact,
$H_1(K;\ZZ)$ is then infinite cyclic and $\aut
(K)$ acts on $H_1(K;\ZZ)$ via $G_K$ with image $\{\pm 1\}$.

\section{Del Pezzo structures on reduced Kodaira curves} 

\subsection*{Homological del Pezzo structures}
Recall that we defined the Lobatchevski lattice $\Lambda_{1,r}$ 
with basis $(l,e_1,\dots ,e_r)$ and $k=3l-e_1-\cdots -e_r$
as a distinguished element. We write $d$ for $9-r$.
If $S$ is a Del Pezzo surface of degree $d$ and $K$ is a reduced  Kodaira curve, then  
the restriction map $\pic (S)\to \pic (K)$ assigns to an 
exceptional curve $E$ on $S$ the class of $E\cdot K$. It has degree one and hence represented by the singleton $E\cap K$, unless $E$ is an irreducible component of $K$. This motivates the following definition.

\begin{definition}\label{homdelpezzo}
Let $K$ be a reduced Kodaira curve. We say that a  homomorphism
$\psi : H_2(K)\to \Lambda_{1,9-d}$  defines a 
\emph{homological Del Pezzo structure} (of degree $d$) on $K$ if
\begin{enumerate}
\item[(i)] the fundamental class $[K]=\sum_{C \text{irr comp}} [C]$  
maps to $k$,
\item[(ii)] if $C,C'$ are distinct irreducible components of 
$K$, then their intersection number in $K$ is equal to 
$\psi [C]\cdot \psi [C']$ and
\item[(iii)] for an exceptional vector $e$ and an irreducible component $C$ of
$K$ we have $e\cdot \psi [C] \ge 0$, unless $\psi [C] =e$, 
\end{enumerate} 
\end{definition}

Notice that if $K$ is irreducible, condition  (ii) is empty and (iii) follows from (i).
The following proposition shows that  this notion helps us to parameterize Del Pezzo triples.

\begin{proposition}\label{prop:chi}
Let  $K$ be a reduced Kodaira curve  and  let $\chi: \Lambda_{1,9-d}\to\pic (K)$ be  a homomorphism such that
\begin{enumerate}
\item[(i)] the dual $\psi: H_2(K)\to \Lambda_{1,9-d}$ of the 
composite $\deg\circ\chi :\Lambda_{1,9-d}\to H^2(K)$  defines a 
homological Del Pezzo structure of degree $d$ on $K$,
\item[(ii)] $\chi$ is nonzero on the roots in that lie in the kernel of 
$\deg\circ\chi$.
\end{enumerate}
Then there is a marked Del Pezzo surface 
$(S,f)$ of degree $d$ and an embedding $K\subset S$ inducing
$\chi$. This system $(S, f,\chi)$ is unique up to unique isomorphism. 
\end{proposition}
\begin{proof}
We first notice that $\chi (l)$ has degree $l\cdot k =3$ on $K$. 
Since $\chi (l)$ has nonnegative degree on every irreducible 
component, the associated linear system is nonempty and maps
$K$ onto a reduced cubic curve in a projective plane
$P$. The image $\bar K\subset P$ is obtained from $K$ by
collapsing each irreducible component on which $l$ has degree zero. 

Suppose an irreducible component $C$ of $K$ has that property. We claim that this
happens  precisely when $C$ is mapped to some $e_i$. For in that case $\psi[C]$ is of the form $\sum_i n_ie_i$. If $\psi[C]\not=e_i$ for some $i$, then  we must have $n_i \le 0$ for all $i$. But the inner product of $\psi[C]$ with an exceptional vector of the form $l-e_i-e_j-e_k$ must also be $\ge 0$. From this we readily see
that all $n_i$ must be zero. This  contradicts the fact that $C$ has nonzero intersection with another irreducible component (in case $K$ is reducible) or represents $K$ 
(in case $K$ is irreducible).
In this situation we  denote  by $p_i\in\bar K$ the point that is the image of $C$ (so with
$\psi[C]=e_i$). The points $p_i\in\bar K$ that we thus obtain are pairwise distinct: if $p_i=p_j$ with $i\not=j$, then $p_i-p_j$ is a root in the kernel of  both $\deg\circ\chi$ and $\chi$.

Any $e_j$ that is not represented by a degree zero component of $K$ defines an element $\chi (e_j)\in \pic (K)$ which has degree $\ge 0$ on each component of $K$. Since the total degree  is $1$,  $\chi (e_j)$ is represented by a unique smooth point of 
$K$. Denote the image of that point in $\bar K$ by $p_j$. The points 
$p_1,\dots ,p_r$ that we now have defined are still pairwise distinct: First of all no two of the newly added points concide since otherwise $\chi$ would vanish on a root in the kernel of $\deg\circ\chi$. And no newly added
$p_j$ equals an old $p_i$ that is associated to an irreducible component $C$ of $K$:
since $e_j.\chi^*(C)=e_j.e_i=0$, $\chi (e_j)$ has zero degree on $C$, hence 
is represented by a point of $K-C$.

Blowing up $p_1,\dots ,p_r$ in $P$ yields 
a surface $S$, a marking  $f:\Lambda_{1,r}\cong\pic (S)$ and an embedding 
$i:K\to S$ such that $\chi=i^*f$. One verifies that $S$ is a Del Pezzo 
surface of degree $d$. The uniquess is left as an exercise.
\end{proof}

Let $K$ be a reduced Kodaira curve,  $p$ a singular point of $K$ and $\psi: H_2(K)\to \Lambda_{1,9-d}$ a homological Del Pezzo structure. This defines in an evident manner  a locus $\Mtilde (K,p;\psi)$ in the moduli space of marked Del Pezzo triples. 

\begin{corollary}
Let $K$ be a reduced Kodaira curve,  $p$ a singular point of $K$ and $\psi: H_2(K)\to \Lambda_{1,9-d}$ a homological Del Pezzo structure. Denote by $Q_\psi\subset \Lambda_{1,9-d}$ the kernel of $\psi^*$. If $K$ is of multiplicative type, then $\Mtilde (K,p;\psi)$ is isomorphic to the toric arrangement complement $\Hom (Q_\psi,\CC^\times)^\circ$ (with the isomorphism given up  to $\pm 1$), unless $K$ is of type $\Ione$,
in which case it is naturally isomorphic to $\Hom (Q_\psi,\CC^\times)^\circ_{\{\pm 1\}}$.
If $K$ is of additive type, then $\Mtilde (K,p;\psi)$ is naturally isomorphic to the projective arrangement complement $\PP(\Hom (Q_\psi,\CC)^\circ$. Both $\Hom (Q_\psi,\CC^\times)^\circ$ and $\PP(\Hom (Q_\psi,\CC))^\circ$ support families of Del Pezzo triples.
\end{corollary}
\begin{proof}
Denote the space of  $\chi\in \Hom (\Lambda_{1,9-d}, \pic (K))$ satisfying the two conditions of Proposition \ref{prop:chi} by $\Hom (\Lambda_{1,9-d}, \pic (K))_\psi^\circ$
so that we have an obvious map $\Hom (\Lambda_{1,9-d}, \pic (K))_\psi^\circ\to \Mtilde (d)$.  This map is essentially the passage to the orbit space of the action of the 
group of  automorphisms of $K$ that  act trivially on $H_2(K)$ and fix $p$.
This group  is easily identified  as $\aut(K)^0$, except when $K$ is of type $\Ione$: it is then all of $\aut (K)$ (which has $\aut(K)^0$ as a subgroup of index $2$).

Denote by $Q_\psi\subset \Lambda_{1,9-d}$ the kernel of $\psi^*$.
We can find for every irreducible component $C$ of $K$ a class $e_C\in \Lambda_{1,9-d}$ such that $\psi^*(e_C)$ takes the value $1$ on $[C]$ and zero on the other irreducible
components. So the $e_C$ define a basis of  $\Lambda_{1,9-d}/Q_\psi$. Since the identity component $\aut(K)^0$ acts simply transitively on $K_\reg$, two elements of
$\Hom (\Lambda_{1,9-d}, \pic (K))_\psi^\circ$  lie in the same $\aut(K)^0$-orbit if and only if they have the same restriction to $Q_\psi$.  Therefore the above map factors
through $\Hom (Q_\psi, \pic(K)^0)^\circ$. The action of $\aut(K)^0$ on $\Hom (Q_\psi, \pic(K)^0)^\circ$ is via  $\pic(K)^0$. This is trivial in the multiplicative case and is
via scalar multiplication in the additive case. All the  assertions of the corollary but the last now follow  easily.  

Proposition \ref{prop:chi} (in its relative form) shows that $\Hom (\Lambda_{1,9-d}, \pic (K))_\psi^\circ$ carries a family of of Del Pezzo triples. Since the group
$\aut(K)^0$ acts freely on the base of this family, the passage to $\aut(K)^0$-orbit spaces still yields a family. This yields the last assertion.
\end{proof}

\subsection*{Glueing constructions} When $K$ is irreducible, both $\psi$ and $p$  become irrelevant and we therefore simply write $\Mtilde (d,\Ione)$ resp.\ $\Mtilde (d,II)$ for the corresponding strata.
We have $\Mtilde (d,\Ione)\cong  \Hom (Q_r, \CC^\times )^\circ_{\{\pm 1\}}=  
(\TT_r^\circ)_{\{\pm 1\}}$ and $\Mtilde (d,II)\cong \PP(\Hom (Q_r,\CC)^\circ$. Notice that in either case the center of $W_r$ acts trivially so that the action is via $\wbar_r$. Both
$\TT_r^\circ$ and  carry $\PP(\Hom (Q_r,\CC)^\circ$ families of Del Pezzo triples.
Recall that the blowup with reflection loci removed, $(\bl_{1}\TT_r)^\circ$, 
comes with an $\sw_r$-action  and is as such the union of the
invariant strata $\TT_r^\circ$ and $\PP(\Hom (Q_r, \CC))^\circ$.

\begin{lemma}\label{lemma:glueII}
The families of Del Pezzo triples over these strata can  be glued to
form one over $(\bl_{1}\TT_r)^\circ$. This family comes with an action of $\sw_r$.
\end{lemma}
\begin{proof}
It is best to start out with a parameterized nodal curve:
the curve $K$ obtained from $\PP^1$ by identifying $0$ and $\infty$, so that
$K_\reg=\CC^\times$. 

Every $e\in\Ecal_r$ defines a homomorphism $\TT_r\to\CC^\times$. 
 If think of $s_e$ as a
section of $K\times \TT_r\to 'TT_r$, then these sections are the restriction of homomorphism from $\Lambda_{1,r}$ to the relative Picard group of $K\times \TT_r /\TT_r$; over $\TT_r^\circ$ this yields a family of degree $d$ Del Pezzo
structures on $K$. But for the moment we rather regard $s_e$ of  $\PP^1\times\TT_r\to\TT_r$. We pull back the latter family with sections to the blowup
$\bl_{1}\TT_r$.  Over the
exceptional locus they all take the value $1\in\CC^\times\subset \PP^1$. Blow up
$\{1\}\times \PP(\Hom (Q_r, \CC))$ in $\PP^1\times \bl_{1}\TT_r$. The
sections $s_e$ will extend across this blow up.
Subsequently we contract the strict transform of 
$\PP^1\times\PP(\Hom (Q_r, \CC))$ in the ambient space to its intersection
with the new exceptional divisor (this is possible, since locally this
this the standard situation over a higher dimensional base).  
Then identify the sections $0$ and $\infty$ (or rather their images) and the result is a family of Kodaira curves over $\bl_{1}\TT_r$.
Its restriction to $\TT_r^\circ$ resp.\ $\PP(\Hom (Q_r, \CC))$ is the
universal family. In the first case this is clear, in the second case
one can verify this easily by restricting to a `linear' curve germ in 
$(\bl_{1}\TT_r)^\circ$ transversal to the exceptional locus.
A relative version of Proposition  \ref{prop:chi} produces a universal family of 
marked Del Pezzo triples over $(\bl_{1}\TT_r)^\circ$, together with
action of $\sw_r$. 
\end{proof}

\subsection*{The type $\Itwo$ case} 
 It is clear that $(\bl_{1}\TT_r)^\circ_{\{\pm 1\}}$ may be identified with 
an open subset of $\tilde M(d)$. We wish to extend this to a description
of $\tilde M(d)$ up to codimension two. This means we will allow also 
the degeneration of a nodal curve into other Kodaira curves. The 
codimension condition  restricts the possibilities to $\Itwo$ and II and 
as we have already dealt with II we now focus on $\Itwo$. 

We begin with determining the homological Del Pezzo structures on such a curve
It is clear that they are permuted by $W_r$ and the
connected component group of $\aut (K)$ and that these actions commute. 
We want to classify the homological Del Pezzo structures 
$\psi : H_2(K)\to \Lambda_{1,9-d}$ with  $d\le 4$ up to this action. 

\begin{lemma} 
Suppose that $K$ is of type $\Itwo$ and let be given a homological 
Del Pezzo structure $\psi :H_2(K)\to \Lambda_{1,9-d}$ of degree $d$. 
Then we have $d\ge 2$ and  there is an element of $W_r$ that  takes the classes of the irreducible components of $K$ to the following systems of vectors in $\Lambda_{1,r}$ (we  also denoted the image of the system in $P_r$):
\begin{enumerate}
\item[(d=2)] $(e_7, k-e_7)\mapsto (\varpi_7, -\varpi_7)$ with $Q_\psi$ of type $E_6$,
\item[(d=3)] $(e_6, k-e_6)\mapsto (\varpi_6, -\varpi_6)$ with $Q_\psi$ of type $D_5$,
\item[(d=4)] $(e_5,k-e_5)\mapsto (\varpi_5, -\varpi_5)$  with $Q_\psi$ of type $D_4$  or\\
$(l-e_1,k-(l-e_1))\mapsto (\varpi_2, -\varpi_2)$ with $Q_\psi$ of type $A_4$,
\item[(d=5)] $(e_4,k-e_4)\mapsto (\varpi_4, -\varpi_4)$  with $Q_\psi$ of type $A_2+A_1$  or\\
$(l-e_1,k-(l-e_1))\mapsto (\varpi_2, -\varpi_2)$ with $Q_\psi$ of type $A_3$,
\end{enumerate}
In brief, the irreducible components of $K$ map to pairs of opposite minuscule weights and each such pair occurs.
\end{lemma}
\begin{proof}
An irreducible component $C$ of $K$ that is not mapped to some $e_i$ is mapped to an element of the form $c=nl-\sum_i n_ie_i$ with $n_0\ge 0$. We also know that 
$c\cdot e_i\ge 0$ for all $i$ and so $n_i\ge 0$. The same must hold
for the other curve $K-C$: $k-c= (3-n)l -\sum_i (1-n_i))e_i$.
So $n\in\{ 0,1,2,3\}$ and $n_i\in\{ 0,1\}$. So there are only
finitely many possibilities for $c$. A straightforward analysis yields the lemma.
\end{proof}

We recall that $\Sigma\subset\Hom(Q_r,\ZZ)$ denotes the set of minuscule weights of 
$Q_r$, that $\TT_r(\Sigma)$ denotes the corresponding toric extension and that
$\hat\TT_r$ stands for $\bl_{\{1\}}(\TT_r)(\Sigma)$. So $\hat\TT_r^\circ$ 
contains  $\bl_{\{1\}}(\TT_r)^\circ$ as an open subset 
and the complement is the union of the $M(\Itwo,\psi)$, where $\Itwo$ represents a fixed Kodaira curve of that type  endowed with a fixed singular point singled out. We have the expected analogue of Lemma \ref{lemma:glueII}:

\begin{lemma}\label{lemma:glueI2}
The families of marked Del Pezzo triples over these strata can  be glued to
form one over $\hat\TT^\circ_r$. This family comes with an action of $\sw_r$.
\end{lemma}
\begin{proof}
As in Lemma \ref{lemma:glueII} we  work with the product $\TT_r\times \PP^1\to \TT^r$
and  the sections $s_e$ ($e$ an exceptional vector) that take values in 
$\TT_r\times \CC^\times$. For every minuscule weight $c$ we have 
one parameter subgroup $\rho_c:\CC^\times 
\to \TT_r$ and a torus embedding $\TT_r\times^{\rho_p}\CC$ which adds to 
$\TT_r$ a copy  $\TT_r(c)$ of the cokernel of $\rho_c$. Given an exceptional vector $e$, then the  corresponding homorphism $s_e:\TT_r\to\CC^\times$ extends to a morphism
$\TT_r\times^{\rho_c}\CC\to \PP^1$ such that its value over the toric 
divisor $\TT_r(c)$ is $0$ if $c(e)>0$,  
lies in $\CC^\times$ if $c(e)=0$, and is $\infty$ if $c(e)<0$.

We now regard each $s_e$ as a section of the trivial $\PP^1$-bundle
over $\TT_r\times^{\rho_p}\CC$. 
Blow up the smooth codimension two locus $\{0\}\times\TT_r(p)$ in 
this bundle and identify  the strict transforms of the constant
sections `zero' and `$\infty$'. The result is a
family of Kodaira curves over $\TT_r\times^{\rho_p}\CC$, of type $\Itwo$
over the toric divisor $\TT_r(c)$ and of type $\Ione$ elsewehere.
The section $s_e$ is now extends to this family when $c(e)\ge 0$: for 
$c(e)= 0$ it will avoid the exceptional divisor and for $c(e)>0$ it
will avoid the strict transform of $\PP^1\times\TT_r(p)$. In either
case, we denote the resulting divisor $D_e$.
\end{proof}

\begin{proof}[Proof of Theorem \ref{thm:moduli}]
Lemma \ref{lemma:glueI2} shows that $\hat\TT^\circ_r$ is the parameter space 
of a family of marked Del Pezzo triples of degree $d=9-r$. The ones that are
missing are those for which the Milnor numbers of
the singular points of $K$ add up to at least $3$. This locus is indeed of 
codimension $\ge 2$. 
\end{proof}

\section{The genus three  case}\label{sect:genus3}

In this section we exclusively deal with the degree two case and so the
relevant root system will be $R_7$ (of type $E_7$). \emph{We  therefore 
often suppress the subscript $7$, it  being understood that $R$
stands for $R_7$, $\TT$ for $\TT_7$ and so on.}

\subsection*{Tacnodal degenerations} We first take up an aspect
on degenerating Kodaira curves that is needed for what follows.

Let $b, c\in \CC [[t]]$ have positive (finite) order 
and consider the double cover of 
$\Kcal\to \PP^1_\Delta$ defined by $y^2=x^2(x^2+2bx+c)$
(here $x$ is the affine coordinate of $\PP^1$). 
We shall assume that $c$ is nonzero and that
$d:=b^2-c$ is not zero. This ensures that the generic
fiber of $\Kcal\to\Delta$ is a nodal curve. 
The fibre $K$ over the closed point is a union of two copies 
$K^\pm$ of $\PP^1$ having in common a single point, where
it displays a {\it tacnodal singularity} 
\index{tacnode} (given by $y^2-x^4=0$).

\begin{definition}\label{def:tacnodal}
We call a  family of reduced Kodaira curves $\Kcal/\Delta$ with general
fiber of type $\Ione$ a \emph{good tacnodal degeneration}\index{tacnodal degeneration, good} if the special fiber is of type II and the family can be given by  the equation as the  double cover of $\PP^1_\Delta$ defined by  $y^2=x^2(x^2+2bx+c)$
with $b, c\in \CC [[t]]-\{ 0\}$ such that $\ord (b^2)>\ord (c)>0$.
\end{definition}

If $d:=b^2-c$ happens to have odd order, then we make the base change defined by
adjoining  $\sqrt{t}$, so that its order becomes even and $\sqrt{d}$ is defined.
We can can then write $x^2+bx+c=(x-u)(x-v)$ with $u,v\in \CC [[t]]$.
Having a good tacnodal degeneration the amounts to $(u+v)^2$ 
having order than $uv$. This comes down to saying that the two roots 
have the same positive order  and have opposite initial coefficients:
$vu^{-1}$ is regular and takes the value $-1$ in $0$. The existence of $\sqrt{d}$ also ensures that the family has two sections $\sigma_{\pm}: t\mapsto (-b\pm\sqrt{d},0)$
whose generic point lies in the smooth part. So we have an isomorphism
$\pic^0(\Kcal_{\Delta^*} )\cong \CC((t))^\times$ (given up to inversion).  

\begin{lemma}\label{lemma:tacnodaldeg}
Suppose we have a good tacnodal degeneration as above that is split is the sense
that $\sqrt{d}\in\CC[[t]]$. 

Then the map from the set of sections of $\Kcal/\Delta$ that meet $K_o$ in
$K_{o,\reg}$ to the
set $\pic_1(\Kcal/\Delta)$ of relative divisor classes of degree one is a bijection
and the specialization map $\pic_1(\Kcal/\Delta)\to \pic_1(K_o)=K_{o,\reg}$ is 
a surjective map of torsors (the underlying surjective  homomorphism
$\pic_0(\Kcal/\Delta)\to \pic_0(K_o)$ is onto).

If $\sigma_1$ and $\sigma_2$ are sections of $\Kcal/\Delta$ meeting $K_o$ in
$K_{o,\reg}$, then the image of $(\sigma_1)-(\sigma_2)$ in $\pic(\Kcal_{\Delta^*} )^0 \cong \CC((t))^\times$ lies in $\pm 1+t\CC[[t]]$ with the sign depending on whether $(+)$ or not $(-)$  $\sigma_1$ and $\sigma_2$ hit the same component of $K_{o,\reg}$.
\end{lemma}
\begin{proof}
The first assertion is well-known and is only included for the record. To prove the second
assertion, 
we introduce the coordinates $\tilde y:=y/x$, $\tilde x:= x+b$, so that
the generic fiber has the equation $\tilde y^2=\tilde x^2-d$, where 
$d=b^2-c$. Notice that the two branches of $\Kcal_{\Delta^*}$  at its singular point have tangents 
$y=x\pm\sqrt{c}$, so that the points $(\tilde x, \tilde y)=(b, \pm\sqrt{c})$ must be omitted.
Since $(\tilde x-\tilde y)(\tilde x+\tilde y)=d$,
$\tilde x+\tilde y$ is a relative paramaterization $\Kcal_{\Delta^*}$.
To be precise, $\tilde x+\tilde y$ maps the smooth part of 
$\Kcal_{\Delta^*}$  to $\PP^1_{\Delta^*}$ with image the complement of $b\pm\sqrt{c}$.
This suggests to use the coordinate
\[
u:= \frac{(\tilde x+\tilde y)-(b-\sqrt{c})}{(\tilde x+\tilde y)-(b+\sqrt{c})}=
\frac{x^2+y+x\sqrt{c}}{x^2+y-x\sqrt{c}}.
\]
instead, for which the two points in question are now $0$ and $\infty$. 
This defines an isomorphism $\pic(\Kcal_{\Delta^*} )^0
\cong \CC((t))^\times$ that maps  $(\sigma_1)-(\sigma_2)$  to  $\sigma_1^*u/\sigma_2^*u$.

The computation is now straightforward. Let $\sigma$ stand for some $\sigma_i$.
We have that  $\sigma(o)$  (assuming it is not a  point at infinity) is of the form 
$(a, \varepsilon a^2)$ with $a\not=0$  and $\varepsilon\in \{\pm 1\}$ (depending on the branch we are on).
So if $\varepsilon=1$, then $\sigma^*u(0)=(2a^2+0)(2a^2-0)^{-1}=1$. If
$\varepsilon=-1$, we need to be a bit more careful. Then 
\[
\sigma^*(x^2+y)s=\sigma^*(x^2-x\sqrt{x^2+2bx+c})=\sigma^*(-bx+c) +\text{higher order terms},
\]
and so 
\[
\sigma^*u=\sigma\frac{(\sqrt{c}-b)\sigma^*x+c}{(-\sqrt{c}-b)\sigma^*+c} +\text{higher order terms}.
\]
Since $b$ has higher order than $\sqrt{c}$, it follows that $\sigma^*u(0)=-1$.

These computations remain valid if $a$ tends to infinity and so the assertion  follows.
\end{proof}

\subsection*{The universal curve of genus three}
Let us begin with a simple observation.
If $C$ is a nonsingular complete genus $3$ curve, then a basis
of its canonical system defines a morphism $C\to \PP^2$.
Denoting the image by $\bar C$, then in
if $C$ is nonhyperelliptic, $C$ maps isomorphically onto $\bar C$ and the image is
a quartic curve. If $C$ is a hyperelliptic, then $\bar C$  is a conic and 
$C\to \bar C$ is dividing out by the hyperelliptic involution. In that case 
the ramification points of $C\to \bar C$ form an $8$-element subset $D$ of $C$.
 If we place ourselves in $\PP^2$, then clearly we can reconstruct $C$ from the pair $(\bar C, D)$ (strictly speaking, up to its 
hyperelliptic involution). It is best regard as $\bar C$ as a conic with multiplicity $2$,
so that its a degenerate quartic. The following lemma is well-known.

\begin{lemma}\label{lemma:thetachar}
Let $\Ccal/\Delta$ be a smooth family genus three curves  whose special  fiber is 
hyperelliptic, but whose generic fiber is not. Let $\Ccal\to \PP^2_\Delta$ be a relative
canonical map and denote  by $\bar\Ccal$  its image and by $D\subset \bar C_o$ the 
ramification set of $C_o\to \bar C_o$. Then every bitangent of the generic fiber of 
$\bar\Ccal/\Delta$  specializes over $o$ to a line through to two distinct points of $D$  and every such  line is thus obtained.
\end{lemma}
\begin{proof}
An odd theta characteristic $\Theta$ of $\Ccal/\Delta$ is smooth of degree $2$ over $\Delta$. The generic fiber $\Theta_{\Delta^*}$ is given by a bitangent of $\bar\Ccal_{\Delta^*}$ in the sense that  $2\Theta_{\Delta^*}$ is the pull-back of a line in $\PP^2_{\Delta^*}$ that is a bitangent  of $\bar\Ccal_{\Delta^*}$ and conversely, any odd theta characteristic of $\Ccal_{\Delta^*}$ so arises. The special fiber $\Theta_o$  is a sum of two Weierstra\ss\ points of $C_o$ and
all pairs of  Weierstra\ss\ points of $C_o$ arise in this manner.
The lemma easily  follows this.
\end{proof}

Let $X\subset \PP^2$ be a nonsingular quartic and let $S\to \PP^2$ be the double cover that ramifies along $X$. Then $S$ is a Del Pezzo surface of degree $2$. The
preimage of a bitangent of $X$ in $S$ consists of two exceptional curves (together forming a Kodaira curve of type $\Itwo$ or III) and we thus get all its $2.28=56$ exceptional curves. Now let $p\in X$ be such that $X$ has a simple tangent at $p$ in the sense
that the projective tangent line $L$ meets $X$ simply in two other points. Then the preimage $\tilde L$ of $L$ in $S$ is an anticanonical curve of type $\Ione$ (its double 
point is $p$). Conversely, any anticanonical curve on $X$ of that type is so obtained.
The 28 bitangents of $S$ meet $L$ in as many points. These are the images of 
points of intersection with $\tilde L$ of exceptional curves of $S$. Since 
these exceptional curves generate $\pic(S)$, we expect that the relative position
of the 28 points on $L$ contains the data that describe the map 
$\pic (S)\to \pic(\tilde L)$. 

In order to understand what happens, let $\Xcal\subset \PP^2_\Delta$  be a family of quartic curves  that is smooth over  $\Delta^*$ and whose closed fiber is a double conic $2C$. We suppose that this family defines a reduced degree 8 divisor $D$ on $C$.
We form the double cover $\Scal\to  \PP^2_\Delta$ with ramification locus $\Xcal$.
Evidently, a  fibers $S_t$ is a degree 2 Del Pezzo surface when $t\not= 0$.
The fiber $S_0$ is the double cover of $\PP^2$ with ramification $2C$ and
hence is the union of two copies (denoted $S_0^\pm$) of $\PP^2$ which meet in $C$ and cross normally along $C$. For every $2$-element subset $A$ of $D$ denote by $E_A\subset \PP^2$ the line passing through $A$. This gives us copies
$E_A^\pm$ in $S_o^\pm$.  It follows from Lemma \ref{lemma:thetachar} that 
each of the 56 elements of the collection $E_A^\pm$ is the specialization of a family $\Ecal_{A}^\pm\subset \Scal$ over $\Delta$ whose generic fiber is an exceptional curve
on the generic fiber of $\Scal/\Delta$.

\begin{lemma}\label{lemma:intersectionlimit}
The intersection numbers in the generic fiber  of  two 
distinct exceptional curves given as above are as follows. Let $A,A'$ be $2$-element subsets of $D$. 
 \begin{itemize}
\item[(i)] If $A\cap A'$ is a singleton, then
$E_{A,t}^\pm\cdot E_{A',t}^\pm =0$ and $E_{A,t}^\pm\cdot E_{A',t}^\mp =1$,
\item[(ii)] if $A\cap A'=\emptyset$, then
$E_{A,t}^\pm\cdot E_{A',t}^\pm =1$ and $E_{A,t}^\pm\cap E_{A',t}^\mp =0$ and
\item[(iii)] $E_{A,t}^+\cdot E_{A,t}^-=2$.
\end{itemize}
\end{lemma}
\begin{proof}
It is a priori clear from Lemma \ref{lemma:thetachar} that such intersection numbers acquire their contributions from
points near $D$. In order to compute  the way a given $p\in D$
contributes, we resolve $\Scal$.
It is not hard to verify that $\Scal$ is smooth outside $D$. To see what happens at 
$p\in D$, we notice that a  local equation for $\Xcal\subset \PP^2_\Delta$ at such a point is $y^2+tx=0$, where $(x,y,t)$ is an analytic coordinate system on a neighborhood of $p$ in $\PP^2_\Delta$. So a local equation for 
$\Scal$ at that point is $y^2+tx=z^2$. This is a cone over a smooth quadric $\Sigma$ in $\PP^3$ and hence gets resolved by a single blowup $\Scalhat\to\Scal$ with a copy of $\Sigma$ as exceptional divisor. The closed fiber
$S_0$ is given by $z=\pm y$. The strict transform  $\hat S^{\pm}_0$ of $S^{\pm}_0$ is a copy of $\bl_p(\PP^2)$ and the three smooth varieties $\hat S^+_0$, $\hat S_0^-$ and $\Sigma$ cross normally.  The intersection $S^\pm_0\cap \Sigma$ corresponds to the exceptional divisor of $\bl_p(\PP^2)\to  S^\pm_0$. If we denote the latter by $F$
and by $o\in F$ the point defined by the  tangent line  of $C$ at $p$, then we may canonically identify $\Sigma$ with $F\times F$ in such a manner that  $S^+_0\cap \Sigma=F\times \{ o\}$
and  $S^-_0\cap \Sigma=\{ o\}\times F$. The involution on $\Scal$ (defined by 
$(x,y,z;t)\mapsto (x,y,-z;t)$) lifts to $\Scalhat$ and acts on $\Sigma=F\times F$
by interchanging the factors.

A pair of exceptional curves $E^\pm_t$ on $S_t$ is the preimage of a bitangent of $X_t$.
If that bitangent specializes to a line $L$ through $p$, then $L$ defines a point 
$[L]\in F-\{ o\}$ and the  exceptional curve $E^+_t$ specializes in $\hat S_0$ to the sum of the strict transform of $E^+_0$ in $\hat S_0^+$  and the rule $[L]\times F$ on 
$\Sigma$. Similarly, $E^-_t$ specializes in $\hat S_0$ to the sum of the strict transform of $E^-_0$ in $\hat S_0^+$  and $F\times [L]$. In particular, the two specializations
meet in our coordinate patch only in $([F],[F])$.  This a point of $\Sigma$ where
$\Scalhat\to \Delta$ is smooth. Since the intersection is simple, the contribution of
$p$ towards $E^+_t\cdot E^-_t$ is $1$. 

For another  pair $E_t'{}^\pm$ of exceptional curves (distinct from $E^\pm_t$) we find
an $[L']\not=[L]$ and we see that the specializations of $E_t^+$ and $E'_t{}^+$ 
in $\hat S_0^+$ do not meet, whereas those of $E_t^-$ and $E'_t{}^+$ simply meet
in $[L']\times [L]$. The same argument as above now shows that the contribution of $p$ towards their intersection number is 
$0$ resp.\ $1$.  This completes the proof.
\end{proof}

\begin{corollary}\label{cor:markinglimit}
A marking of $\Scal_{\Delta^*}$ is obtained by numbering the points of $D$:
$D=\{p_0,\dots ,p_7\}$ and the rule
\begin{align*}
e_i&\mapsto E_{p_0,p_i}^-, \quad 1\le i\le 7\\
(2l-e_1-\cdots -e_7)+e_i+e_j&\mapsto 
E_{p_i,p_j}^-,\quad 1\le i<j\le 7,\\
l-e_i-e_j&\mapsto E_{p_i,p_j}^+,\quad 1\le i<j\le 7,\\
(3l-e_1-\cdots -e_7)-e_i&\mapsto 
E_{p_0,p_i}^+,\quad 1\le i\le 7.
\end{align*}
\end{corollary}
\begin{proof}
Following Lemma \ref{lemma:intersectionlimit} this assignment respect the inner products. Since the collection of exceptional vectors generates $\Lambda_{1,7}$, the  
first part of the corollary follows. 
\end{proof}

Consider the map that sends $E^\pm_A$ to $\pm 1$. Via the marking
defined by corollary \ref{cor:markinglimit} we see that it is simply given by the 
parity of the coefficient of $l$. It therefore defines a
homomorphism $\Lambda _{1,7}\to \{\pm1\}$.
Now observe that the restriction of that homomorphism to $Q_{}$ (which is an element 
of $\Hom(Q_{},{\pm1})=\Hom(Q_{},\CC^\times)\subset \TT_{}$)  is the 
weight $e^{\varpi_1}$ defined by the nonspecial end vertex of the $\hat E_7$-diagram.
Its stabilizer in $W_{}$ is $\{ \pm1\} .W_{\varpi_1}$, where $W_{\varpi_1}$ is of type $A_{}$ and hence isomorphic tot the symmetric group of 8 elements.
This group has two orbits in $\Ecal_7$ that are distinguished by the parity of the coefficient of $l$ and in terms of the indexing by the  $E^\pm_A$'s by the sign.
\\

Let us now fix $p\in C$ distinct from the support of the associated degree 8 divisor $D$ on $C$ and let  $\gamma$ be a section of $\Xcal/\Delta$  with $\gamma (o)=p$.

\begin{lemma}\label{lemma:gooddeg}
The generic fiber  $X_{\Delta^*}$ has a simple tangent at $\gamma({\Delta^*})$ and if we denote that line by $\Lcal_{\Delta^*}\subset \PP^2_{\Delta^*}$, then  this family extends to 
$\Lcal/\Delta$ with the fiber over $o$ being the tangent line of $C$ at $p$.
Moreover, the preimage  $\Kcal/\Delta\subset \Scal/\Delta$ of $\Lcal/\Delta$
defines a good tacnodal degeneration in the  sense of
Definition \ref{def:tacnodal}.
\end{lemma}
\begin{proof}
We choose affine coordinates $(x,y)$ for $\PP^2_\Delta$ such 
that $p=(0,0;0)$ and $C$ is given by $(y-x^2)^2=t=0$. So
$\Xcal$ has an equation of the form $F(x,y)=0$, with 
\[
F(x,y)=(y-x^2)^2+t^r G(x,y), \quad  r\ge 1,
\]
where $G\in \CC[[t]][x,y]$ is a quartic polynomial in $(x,y)$ that is not divisible by $t$ and 
$G_o(x,y)$ is not divisible by $y-x^2$.
The degree 8 divisor on $C$ is defined by $G_0= 0$ and so we must have 
$G_o(0,0)\not= 0$. So if $\gamma=(\gamma_x ,\gamma_y)$, with $\gamma_x,\gamma_y\in \CC[[t]]$, then the identity
\[
(\gamma_y-\gamma_x^2)^2+t^r G(\gamma_x,\gamma_y)=0
\]
shows that that $r$ must be even and that 
$\gamma_y-\gamma_x^2$ has order $r/2$. The line through $(0,0)$ that is parallel to the tangent line of $X_t$ at $\gamma $ is given by 
\[
2(\gamma_y-\gamma_x^2)(dy-2\gamma_x)dx) +t^r dG(\gamma_x,\gamma_y)=0
\]
If we divide this expression by $t^{r/2}$ and the put $t=0$, we find the line with equation
$dy=0$. This is the slope of the limiting line. That line also passes through $(0,0)$ and so the first part of the lemma follows.

For the second part, we parameterize $\Lcal$ in the obvious way:
$u\mapsto (\gamma_x  +u, \gamma_y +\alpha(t) u)$ where we know that $\alpha $ has order $\ge 1$. 
Substitution of this parameterization in $F(x,y;t)$ gives
\begin{multline*}
F(\gamma_x +u, \gamma_y+\alpha u)=
((\gamma_y +\alpha u)-(\gamma_x +u)^2)^2+t^r G(\gamma_x+u, \gamma_y+\alpha u).
\end{multline*}
The lefthand side must be divisible by $u^2$ and will therefore have the form
$u^2(a u^2+b u+c)$ with $a,b,c\in \CC((t))$.  Modulo $t^r$
this must be equal to
\[
\Big(\gamma_y +\alpha u)-(\gamma_x +u)^2\Big)^2 =
\Big((\gamma_y -\gamma_x^2) +(\alpha -2\gamma_x )u -u^2\Big)^2
\]
and thus we find that 
\begin{align*}
a & \equiv -1 \pmod{t^r}, \\
b & \equiv -2(\alpha (t)-2\gamma_x )  \pmod{t^r}, \\
c & \equiv -2(\gamma_y -\gamma_x  ^2)+(\alpha -2\gamma_x )^2 \pmod{t^r}, \\
0 & \equiv 2(\alpha -2\gamma_x )(\gamma_y -\gamma_x ^2) \pmod{t^r}.
\end{align*}
So $a$ has order zero.
Since $\gamma_y -\gamma_x ^2$ has exact order $r/2$, it follows from the 
last equation that $\alpha-2\gamma_x $ has order $> r/2$.
Hence $b/2a$ has order $> r/2$ and $c/a$ has exact order $r/2$
and so the last assertion of the lemma  follows.
\end{proof}

We can now relate this to our torus $\TT_{}$. 

\begin{lemma}
If $A,A'$ are $2$-element subsets of $D$ whose intersection is a singleton, then
Lemma \ref{lemma:intersectionlimit} associates to $E^+_A-E^+_{A'}$ a root (denoted $\beta (A,A')$) whose reflection fixes $e^{\varpi_1}$. All roots with that property are thus obtained and make up a subsystem of type $A_7$ that has 
\[
(\beta_1,\dots ,\beta_7):=(e_2+\cdots +e_7-2l, e_1-e_2,e_2-e_3,\dots ,e_6-e_7)
\] 
as a root basis.
\end{lemma}
\begin{proof} With the help of Corollary \ref{cor:markinglimit} we easily verify that
$\beta (A,A')$ runs over all the roots for which the coefficient of $l$ is even.
These are precisely the roots  whose reflection fixes $e^{\varpi_1}$. The rest 
is straightforward.
\end{proof}

Recall that our marking determines a group homomorphism $\chi :Q_7\to \pic(\Kcal_{\Delta^*})^0\cong \CC((t))^\times$.

\begin{proposition}\label{prop:a7limit}
For a section $\gamma$ of $\Xcal/\Delta$ as above, the
resulting family of Del Pezzo triples over ${\Delta^*}$ defines an ${\Delta^*}$-valued
point of $\TT_{}^\circ$
that extends to $\Delta\to \TT_{}$ with $o$ mapping to  $e^{\varpi_1}$.  
Its lift  to $\bl_{e^{\varpi_1}}(\TT_{})$ sends $o$ to a point  that does not lie on the strict transform of a reflection hypertorus (hence is in $\tilde\TT^\circ_{}$) and
which describes the relative position of $p_0,\dots ,p_7$ on the punctured conic
$C-\{p\}$ (which we think of as an affine line) as follows:  the mutual ratio of $[p_0-p_1:\cdots :p_6-p_7]\in\PP^7$ is the value  of $[\chi (\beta_1)-1:\cdots ;\chi (\beta_7)-1]$ in
$o$. 
\end{proposition}
\begin{proof}
We have seen that each $E^\pm_A$ is the specialization of a unique exceptional curve 
$\Ecal^\pm_A$ over $\Delta$. The latter meets $\Kcal$ in a section of $\Kcal/\Delta$
that we shall denote by  $\sigma^\pm_A$. The value of this section in $o$ lies over the point $p_A$ where the line through $A$ and the tangent to $C$ at $p$ meet and since $p_A\not= p$, this value is a smooth point of $K_o$.
Two such sections hit the same  same component of $K_{o,\reg}$ if and only they
have the same sign in their labeling and in that case, according to  Lemma 
\ref{lemma:tacnodaldeg}, their difference in  $\pic(\Kcal/\Delta)$ then specializes over $o$ to $1$ (it will be $-1$ otherwise). This implies the first part of the proposition:  the closed point of $\Delta$ will go to a point of $\TT_{}$ where all the characters defined by exceptional vectors take only two (opposite) values, the value only depending 
on  the parity of the coefficient of $l$ in the exceptional vector.  Since roots are differences of exceptional vectors, this means that
the character defined by any root takes there the  value $\pm 1$, 
depending on the parity of $l$ in the root. So this must be $e^{\varpi_1}$. 

If $A\cap A'$ is a singleton, then the lines through $A$ and $A'$ clearly do not meet outside $A\cap A'$ and hence $p_A\not=p_{A'}$. Our marking associates to such a 
pair the root $\beta_{A,A'}\in R_{}$.
If we fix $A'$ for the moment, then  the entries of $(p_A-p_{A'})_{|A\cap A'|=1}$ are all nonzero (these differences makes sense if we regard  $L_o-{p}$ as an affine line)
and it follows from the first part of \ref{lemma:tacnodaldeg} that their mutual ratio
is equal to that  the collection $\big((\sigma^+_A)-(\sigma^+_{A'})-1\big)_{|A\cap A_0|=1}$ evaluated in $o$. But then this  remains true if we let both $A$ and $A'$ vary subject to the condition $|A\cap A'|=1$. For such a pair $(A,A')$ we have
$(\sigma^+_A)-(\sigma^+_{A'})=\chi(\beta_{A,A'})$ and since any root whose reflection
fixes $\varpi_0$ is of the 
form $\beta_{A,A'}$ all assertions of the proposition but the last follow.  

To prove the last assertion, we note that assigning to a point 
$a\in C-\{p\}$ the point $a'$ of intersection of its tangent with $L_o$ defines an isomorphism $C-\{p\}\cong L_o-\{p\}$. For $a,b\in C-\{ p\}$, the line spanned by 
$a,b$ (the tangent to $C$ if $a=b$) meets $L-\{ p\}$ in 
a point which is easily verified to be the barycenter of $a'$ and $b'$. Thus 
$E_{p_i,p_j}\cap L_o =\frac{1}{2}(p'_i+p'_j)$. Now $E_{p_i,p_j}\cap L_o$ is the
image of $E^+_{p_i,p_j}\cap K_o=\sigma_{p_i,p_j}^+(o)$ and this is associated to the exceptional vector $e_j$ if $0=i<j$ and to $(2l-e_1-\cdots -e_7)+e_i+e_j$ if $1\le i<j$.
This makes $\beta_1=e_2+\cdots +e_7-2l=e_2-((2l-e_1-\cdots -e_7)+e_1+e_2)$ correspond to $\frac{1}{2}(p'_0+p'_2)-\frac{1}{2}(p'_1+p'_2)=\frac{1}{2}(p'_0-p'_1)$ and
$\beta_i=e_i-e_{i+1}$ for $i=1,\dots ,6$ to $\frac{1}{2}(p'_0+p'_i)-\frac{1}{2}(p'_0+p'_{i+1})=\frac{1}{2}(p'_i-p'_{i+1})$. The last clause follows.
\end{proof}

Recall that we have a family of marked Del Pezzo triples of $\hat\TT^\circ$ with the
underlying family of  Del Pezzo surfaces being a
double cover over $\PP^2_{\hat\TT^\circ}$. The ramification locus 
$\Xcal_{\hat\TT^\circ}\subset \PP^2_{\hat\TT^\circ}$ of that double cover is then a family of  nonsingular  quartic curves over  $\hat\TT^\circ$. It comes with a section, making it
a family of pointed quartic curves. We can now fill in a family of pointed hyperelliptic genus $3$  curves over $\tilde\TT^\circ -\hat\TT^\circ$ as to obtain a family of pointed smooth genus $3$ curves over $\tilde\TT^\circ$.

\begin{corollary}\label{cor:hyperellipticfill}
The closure of $\Xcal_{\hat\TT^\circ}$ in $\PP^2_{\tilde\TT^\circ}$ yields after normalization
a family of pointed genus $3$ curves $\Xcal_{\tilde\TT^\circ}$ over $\tilde\TT^\circ$, which is over  $\tilde\TT^\circ-\hat\TT^\circ$ hyperelliptic. The section extends across   $\tilde\TT^\circ$ and avoids the  Weierstra\ss\ points.
\end{corollary}
\begin{proof}
It suffices to deal with the stratum over $e^{\varpi_1}$.  This stratum can be identified with
$\PP(V^\circ)$, where $V$ is the span of an $A_7$ root system with root basis 
$\beta_1,\dots ,\beta_7$. An element of $\PP(V^\circ)$ defines an affine line 
with a configuration of $8$ distinct numbered points $p_0,\dots ,p_7$, given up to translation. We take as affine line a punctured conic $C-\{ p\}\subset \PP^2$.
If $F\in \CC[Z_0,Z_1,Z_2]$ is an equation for $C$ and choose $G\in \CC[Z_0,Z_1,Z_2]$
homogeneous of degree $4$ such that $C\cap (G=0)=\{p_0,\dots ,p_7\}$. Then
$F^2+tG$ defines a degenerating  family of quartic curves  $\Xcal/\Delta$ whose normalization $\tilde\Xcal/\Delta$ is a family of genus three curves. Choose also a section 
$\gamma$ of  $\tilde\Xcal/\Delta$ such that $\gamma(o)=p$. According to Proposition \ref{prop:a7limit} these data define  a morphism $\Delta\to \tilde\TT^\circ$ whose generic point lies in $\hat\TT^\circ$ and whose special point goes to the given element of $\PP(V^\circ)$. The universal nature of this construction yields the corollary.
\end{proof}

\begin{proof}[Proof of Theorem \ref{thm:moduligenus3}] 
Most of what the theorem states is immediate from Corollary \ref{cor:hyperellipticfill}.
The assertion that remains to prove is that the complement of
$(\tilde\TT^\circ_{})_{SW_{}}$ in $\Mcal_{3,1}$ is of  codimension $\ge 2$. That complement parameterizes pointed genus three curves $(X,p)$  for which, 
in case is $X$ a quartic curve, $p$ is a hyperflex point and in case $X$ is hyperelliptic, 
$p$ is a Weierstra\ss\ point (in either case this amounts to: $4(p)$ is a canonical divisor).
It it clear that this defines a locus in  $\Mcal_{3,1}$ of codimension $\ge 2$. 
\end{proof}

\begin{remark}
If $S$ is a Del Pezzo surface of degree $2$ and $X\subset S$ is its fixed point set of the involution, then a marking $\Lambda_{1,7}\to \pic (S)$ of $S$ composed with the restriction map yields a homomorphism $\Lambda_{1,7}\to \pic (X)$. Now an
odd theta characteristic of $X$ is of the form $(p)+(q)$, with $2(p)+2(q)$ a canonical divisor and so this is where $X$ meets a bitangent (relative to its canonical embedding).
Such a bitangent corresponds to a pair of exceptional curves on $S$ whose
sum is anticanonical  and hence this
homomorphism sends an exceptional vector to an odd theta characteristic of $X$ and 
$k$ to the canonical class $[K_X]\in \pic_4(X)$.  We thus obtain a bijection between the set $\Theta_\odd(X)\subset\pic^2(X)$ of odd theta characteristics of $X$ and pairs
of exceptional vectors with sum $k$. 

An exceptional vector defines an element of the weight lattice $Q^*_{}:=\Hom(Q_{},\ZZ)$ and the bijection is in fact the restriction of an map between two torsors: 
namely from $Q_{}[1]:=\{ v\in\Lambda_{1,7}\, |\, v\cdot k=1\}$ (a $Q_{}$-torsor) to
the set $\Theta (C)$ of all theta  characteristics of $X$ (a torsor under the group of order two elements in $\pic(X)$, $\pic(X)_2$). This map is surjective and factors through an isomorphism of torsors
\[
Q_{}[1]/2Q^*\cong \Theta (C).
\]
The underlying group isomorphism $Q/2Q^*\cong \pic(X)_2$ had been noted by Van Geemen as well as the fact that this is an isomorphism of symplectic modules:
the evident $\ZZ /2$-valued bilinear form on $Q/2Q^*$ corresponds to the natural pairing on $\pic(X)_2$. This also identifies  $\overline{W}=W_{}/\{ \pm 1\}$
with the symplectic group of $\pic(X)_2$. So the variety 
$\tilde\TT^\circ_{\{\pm 1\}}$ can be regarded is an open subset in the moduli space 
$\Mcal_{3,1}[2]$ of pointed genus three curves with level two structure.

The marked Del Pezzo surfaces of degree $2$ define  a connected degree two covering
of a connected component of $\Mcal_{1,3}[2]$ and a $W$-covering of $\Mcal_{3,1}$
(which has $\tilde\TT^\circ$ as  an open set). This is not an unramified covering, even if $\Mcal_{1,3}$ is regarded as a stack, but one may wonder whether this still has some meaning in terms of  a structure on $X$.
\end{remark}

\subsection*{The genus three surface with a boundary circle}
Recall that $\Mcalvec$ denotes the moduli space of pairs $(X,v)$, where $X$ is projective nonsingular genus three curve and $v$ a nonzero tangent vector of $X$
and that it comes with a proper $\CC^\times$-action whose orbit space is $\Mcal_{3,1}$. 
We have indentified $\tilde\TT^\circ_W$ with an open subset of
$\Mcal_{3,1}$ whose complement has  codimension $\ge 2$. So it enough to determine
the fundamental group of the restriction of $\Mcalvec$ to $\tilde\TT^\circ_W$.
We find it convenient to work $W_{}$-equivariantly with the pull-back of $\Mcalvec$ to
$\tilde\TT^\circ_{}$. Over that variety  it is  a genuine $\CC^\times$-bundle  with $W_{}$-action that we shall denote by $\Vcal^\times \to \tilde\TT^\circ_{}$; we reserve the symbol
$\Vcal$ for the associated line bundle over $\tilde\TT^\circ_{}$.

Let us begin with the following observation. If $(S,K,p)$ is a Del Pezzo triple of degree
two, then the natural involution $\iota$ of $S$ has the property that its fixed point set  $X\subset S$ is a quartic curve passing through $p$ so that $T_pX$ is the
fixed point line of $\iota$ in $T_pS$. Since $T_pS$ is also the Zariski tangent space 
$T_pK$ of $K$ at $p$, we can reconstruct $T_pX$ from the triple $(X,p,\iota )$.
Our  construction of $\TT^\circ$   yields a trivial bundle of type $\Ione$ curves with involution over $\TT^\circ$: if $K$ denotes such a curve with involution $\iota$ so
that the family is $K\times \TT^\circ\to\TT^\circ$,  then the action of  $W_{}$ respects
this decomposition with $w\in W$ acting on $K$  by $\iota$ if $\det (w)=-1$.
So from the above discussion we see that $\Vcal \big|_{\TT^\circ}$ is simply
$(T_pK)^\iota\times \TT^\circ \to \TT^\circ_{}$ with $W$ acting trivially on the first factor. In particular,  $(\Vcal^\times\big| \TT^\circ)_W$ can be identified with 
$\CC^\times \times \TT^\circ_W$. Since $(\Vcal^\times\big| \TT^\circ)_W$ is open-dense
in $\Vcal^\times$, the fundamental group of $\Vcal^\times_W$ (which is also the one of 
$\Mcalvec$) is a quotient of the direct product of an infinite cyclic group $r^\ZZ$ (here we think of $r$ as the positive generator of $\pi_1(\CC^\times,1)$) and the orbifold fundamental group $\{1,\iota_1\} \ltimes \Ar_{\hat E_7}$ of $\TT^\circ_W$ (here $\iota_1$ denotes
the nontrivial involution of the $\hat E_7$-diagram). In fact, the complement of $(\Vcal^\times\big| \TT^\circ)_W$ in $\Vcal^\times_W$ consist of three connected
divisors, and each of these defines a relation in $r^\ZZ\times (\{1,\iota_1\} \ltimes \Ar_{\hat E_7})$ so that  $\pi_1(\Vcal^\times)$ is obtained by imposing these relations. So we need to identify these relations.

To this end, we write the  complement $\tilde\TT^\circ_{}-\TT^\circ_{}$ as a sum of three of divisors:
$D(E_7)$,  the divisor over the unit element of $\TT_{}$  (and over which $K$ degenerates into a type II curve), $D(A_7)$ the sum of the divisors over the $A_7$-points (over which $K$ degenerates into a type III curve) and $D(E_6^\tor)$, the sum of the toric divisors (over which $K$ degenerates into a type $\text{I}_2$ curve).

\begin{lemma}
The $W$-line bundle $\Vcal$  over  $\tilde\TT^\circ$ is $W_{}$-equivariantly 
isomorphic to 
\[
\Ocal_{\tilde\TT^\circ}(-2D(E_7)-D(A_7)+D(E_6^\tor)).
\]
\end{lemma}
\begin{proof}
Since $W_{}$ acts transitively on the connected components of a divisor of a given type,
the bundle $\Vcal$ will a priori be given by a linear combination of the three divisors.
Our task is therefore to find the coefficients. We do this with a valuative test:
we pull back along a morphism $\gamma:\Delta\to \tilde\TT^\circ_W$ whose special point goes to one of the added divisors and which is  transversal to that divisor. Since $\tilde\TT^\circ\to \tilde\TT^\circ_W$ has order two
ramification along these divisors (they are pointwise fixed under the involution $\iota$),  
a base change of order two ($\tilde t^2=t$) yields a lift $\tilde\gamma :\tilde \Delta\to  \tilde\TT^\circ$.
The resulting degeneration  $\Kcal/\tilde\Delta$ is split ($\Kcal_{{\tilde\Delta^*}/{\tilde\Delta^*}}$ is trivial
in the sense that it is ${\tilde\Delta^*}$-isomorphic to $K\times_{\text{Spec}(\CC)}{\tilde\Delta^*}$)  and an affine part of $\Kcal$ admits the following weighted homogeneous  equation:
\begin{enumerate}
\item[(II)] $y^2=x^3+\tilde t^2x^2$ with $(x,y,\tilde t)$ having weights $(2,3,1)$, ($\omega=\tilde t^{-2}dx$);
\item[(III)] $y^2=x^4+\tilde t^2x^2$ with $(x,y,s)$ having weights $(1,2,1)$,  ($\omega=\tilde t^{-1}dx$);
\item[($\text{I}_2$)] $y^2=x^2+\tilde t^2x^4$ with $(x,y,t)$ having weights $(1,1,-1)$,
 ($\omega=\tilde tdx$).
\end{enumerate}
Here the form $\omega$ has (and is up to a scalar in $\CC^\times$ characterized by) the following three properties: it is $\iota$-invariant, it is nonzero on the  Zariski tangent space of the singular point of $\Kcal_{{\tilde\Delta^*}}$ and  it has weight zero. The last property ensures that if we trivialize $K_{{\tilde\Delta^*}/{\tilde\Delta^*}}$, then $\omega$ gets trivialized, too. So $\omega$
might be regarded as the pull-back of a constant section of $\Vcal^*\big|_{\TT^\circ}$.
On the other hand in all three  cases the form $dx$ is $\iota$-invariant and nonzero on the  Zariski tangent space of the singular point of the central fiber as well.
So in the first case, $\tilde t^2\omega$ defines an isomorphism of the pull-back
of $\Vcal$ with $\tilde t^2\CC[[\tilde t]]$ and likewise for the other cases. The lemma follows.
\end{proof}

We will find it convenient to work on the quotient of $\Vcal^\times$ by the involution
given by scalar mulitplication in the fibers by $-1$. This quotient is of course 
simply $(\Vcal\otimes\Vcal)^\times$. In this way $r^\half$ can be understood
as an element of the orbifold fundamental group of $(\Vcal\otimes\Vcal)^\times_W$
and we identify the fundamental group of  $\Vcal^\times_W$ with a subgroup of index two of the orbifold fundamental group of $(\Vcal\otimes\Vcal)^\times_W$.

\begin{lemma}\label{lemma:genus3relations}
The orbifold fundamental group of $(\Vcal\otimes\Vcal)^\times_W$ is obtained from
$r^{\half\ZZ}\times \sar_{\hat E_7}$ (the fundamental group of  
$(\Vcal\otimes\Vcal|_{\TT^\circ})^\times_W$) by imposing the following three relations
(corresponding to the three added divisors):
\begin{align*}
D(E_7):&\;  (r,\Delta_{E_7})\equiv 1;\\
D(A_7):&\; (r^\half, \iota_1\Delta_{A_7})\equiv 1;\\
D(E_6^\tor):&\; (r^{-\half},(\iota_1 \Delta_{E_6})^{-1})\equiv 1.
\end{align*}
\end{lemma}
\begin{proof}
We use our fixed trivialization  over of $\tilde\TT^\circ$.

Let us begin with  the first case, where the divisor is $D(E_7)$.
In terms of our fixed trivialization, the
section $\tilde t=1$ of $\tilde\gamma^*\Vcal^\times$ is over $\tilde\Delta^*$ given by $(\tilde t^2,\tilde t)\in \CC^\times\times\tilde\Delta^*$. The homotopy class 
of this restriction is $(r^2,\Delta_{E_7}^2)$ and so this represents a trivial element
of the fundamental group of $\tilde\TT^\circ$. But $\tilde\gamma^*\Vcal^\times$
is obtained from $\gamma^*\Vcal^\times_W$ by the substitution $\tilde t^2=t$ and so 
$(r,\Delta_{E_7})$ represents a trivial element of the
fundamental group of $\Vcal^\times_W$.

The second case is dealt with in the same way, except that
we now work in $r^{\half\ZZ}\times \sar_{\hat E_7}$ (the orbifold fundamental group of $(\Vcal\otimes\Vcal)^\times_W$).

The difference between the third case and the previous two cases is that
we need to explain why we get in the second factor not $\iota_1\Delta_{E_6}$,
but its inverse. The answer to that is that in $\CC^\times$ simple positive loops around $0$ and $\infty$ represent in the fundamental group of $\CC^\times$ elements that are each others inverse. To make this precise, consider the weight 
$\varpi_7\in\Hom (Q,\ZZ)$ (given by taking  the inner product with $e_7$). Then
we have a cocharacter $e^{\varpi_7}: \CC^\times \to \Hom (Q,\CC^\times)=\TT$.
The preimage of $\TT^\circ$ is $\CC^\times-\{ 1\}$ and the composite
map $\CC^\times \to\TT\to\TT_W$ factors through a closed embedding
of the quotient of $\CC^\times$ by the involution $z\mapsto z^{-1}$.  Let us denote the orbit space of $\PP^1$ by this involution $C$ and denote the image of $0$ resp. $1$ by
$p_0$ resp.\ $p_1$. Then a  simple loop around $0$ maps to simple loop around $p_0$, whereas a simple loop around $1$ maps to the square a simple loop around
$p_1$. Since a simple positive loop around $p_1$ represents $\iota_1\Delta_{E_7}$, a simple positive loop  around $0$ represents its inverse $(\iota_1\Delta_{E_7})^{-1}$.
This completes the proof.
\end{proof}

\begin{proof}[Proof of Theorem \ref{thm:genus3}]
The fundamental group of $\Mcalvec$ is that of $\Vcal^\times_W$ and Lemma 
\ref{lemma:genus3relations} shows that this equals the image of $r^\ZZ\times \sar_{\hat E_7}$ in the quotient of  $r^{\half\ZZ}\times \sar_{\hat E_7}$ by the three relations of that lemma. It is clear that this is simply the quotient $\sar_{\hat E_7}$ by the relations
$\Delta_{A_7}\equiv\Delta_{E_6}$ and $\Delta_{E_6}^2\equiv\Delta_{E_7}$.
\end{proof}

\end{document}